\documentclass[11pt, a4paper]{amsart}
\usepackage[latin2]{inputenc}
\usepackage[T1]{fontenc}
\usepackage{lmodern}
\usepackage{amsmath, amsthm,amssymb}
\usepackage{array,tabularx}
\usepackage{enumerate}
\usepackage{pict2e}
\usepackage[
      bookmarksnumbered, bookmarksopen,
      naturalnames,
]{hyperref}

\swapnumbers
\theoremstyle{plain}
\newtheorem{theorem}{Theorem}[subsection]
\newtheorem{lemma}[theorem]{Lemma}
\newtheorem{corollary}[theorem]{Corollary}
\newtheorem{proposition}[theorem]{Proposition}
\newtheorem{cordef}[theorem]{Corollary/Definition}

\theoremstyle{definition}
\newtheorem{definition}[theorem]{Definition}

\newtheorem{disdef}[theorem]{Discussion/Definition}
\theoremstyle{remark}
\newtheorem{remark}[theorem]{Remark}
\newtheorem{example}[theorem]{Example}
\newtheorem{terminology}[theorem]{Terminology}
\newtheorem*{fact}{Fact}
\newtheorem{termdis}[theorem]{Terminology/Discussion}
\newtheorem{notation}[theorem]{Notation}

\swapnumbers
\newtheoremstyle{para}{}{}{\normalfont}{}{}{}{0pt}{\thmnumber{#2}.\thmname{#1} \textbf{\thmnote{#3. }}}
\theoremstyle{para}
\newtheorem{para}[theorem]{}

\newenvironment{fakeequation}{\begin{equation}\begin{minipage}{0.7\linewidth}}{\end{minipage}\end{equation}}

\numberwithin{equation}{section}

\def\MRgetnumber#1 #2\relax{\href{http://www.ams.org/mathscinet-getitem?mr=#1}}
\newcommand*{\itemtitle}[1]{\textbf{#1}}

\newcommand*{\combinatorialarea}[1]{\ensuremath{\mathinner{g(#1)}}}
\newcommand*{\facevalue}[1]{\ensuremath{\mathinner{m\sb{#1}}}}

\newcommand*{\selfintersectionnumber}[1]{\ensuremath{\mathinner{b\sb{#1}}}}
\newcommand*{\Eulernumber}[1]{\ensuremath{\mathinner{e\sb{#1}}}}
\newcommand*{\innerdeterminant}[2]{\ensuremath{\mathinner{n\sb{#1,#2}}}}
\newcommand*{\subchaindeterminant}[2]{\ensuremath{\mathinner{c\sb{#1,#2}}}}
\newcommand*{\multiplicity}[2]{\ensuremath{\mathinner{t\sb{#1,#2}}}}

\newcommand*{\legdeterminant}[2][]{\ensuremath{\mathinner{n\sb{#2}\sp{#1}}}}

\newcommand*{\namedvector}[1]{\ensuremath{\mathinner{\overrightarrow{\mathbf{#1}}}}}

\newcommand*{\normalvector}[1]{\ensuremath{\mathinner{{\namedvector{a}}\sb{#1}}}}
\newcommand*{\coordinatevector}[1]{\namedvector{e\sb{#1}}}

\newcommand*{\vectorbycoordinates}[1]{\ensuremath{\mathinner{(#1)}}}

\newcommand*{\scalarproduct}[2]{\ensuremath{\mathinner{\left\langle
      #1, #2
    \right\rangle}}}
\newcommand*{\coordinateof}[2]{\ensuremath{\mathinner{\scalarproduct{#2}{\coordinatevector{#1}}}}}

\newcommand*{\trianglebyvertices}[3]{\ensuremath{\mathinner{\triangle\sb{#1 #2 #3}}}}

\newcommand{\mygcd}[1]{\ensuremath{\mathinner{\gcd\left(#1\right)}}}

\newcommand{\setZ}{\mathbb{Z}}
\newcommand{\setQ}{\mathbb{Q}}
\newcommand{\setC}{\mathbb{C}}
\newcommand{\setR}{\mathbb{R}}
\newcommand{\setN}{\mathbb{N}}
\newcommand{\setT}{\mathbb{T}}

\newcommand{\calS}{\mathcal{S}}
\newcommand{\calD}{\mathcal{D}}
\newcommand{\calf}{\mathcal{F}}
\newcommand{\calr}{\mathcal{N}}
\newcommand{\calv}{\mathcal{V}}

\DeclareMathOperator*{\supp}{supp}

\providecommand{\coloneqq}{\mathrel{:=}}
\newcommand{\ep}{\ensuremath{{\mathrm{M1}}\sb{+}}}
\newcommand{\kp}{\ensuremath{{\mathrm{M2}}\sb{+}}}
\newcommand{\en}{\ensuremath{{\mathrm{M1}}\sb{-}}}
\newcommand{\kn}{\ensuremath{{\mathrm{M2}}\sb{-}}}
\newcommand{\jn}{\ensuremath{{\mathrm{Mj}}\sb{-}}}
\newcommand{\epn}{\ensuremath{{\mathrm{M1}}\sb{\pm}}}
\newcommand{\kpn}{\ensuremath{{\mathrm{M2}}\sb{\pm}}}
\newcommand{\jpn}{\ensuremath{{\mathrm{Mj}}\sb{\pm}}}

\newcommand{\move}[1]{\xrightarrow{#1}}\newcommand{\nrt}{\texorpdfstring{\blacksquare}{Trapezoid}}
\newcommand{\ctr}{\texorpdfstring{\blacktriangle}{Triangle}}
\newcommand{\ced}{\texorpdfstring{{\mathbf{l}}}{Edge}}

\title{Invariants of Newton
  non-degenerate surface singularities}
\keywords{hypersurface singularities, links of singularities,
  resolution graphs, Newton boundary, Newton polyhedrons.}
\subjclass[2000]{Primary: 14J17, 14Q10;
  Secondary: 52B20}
\author{Gábor Braun}
\email{braung@renyi.hu}
\urladdr{\url{http://www.renyi.hu/~braung}}
\thanks{The first author is partially supported by Hungarian National
    Research Fund, grant No. T 042 769.}
\address{Alfréd Rényi Institute of Mathematics\\
    1053 Budapest,  Reáltanoda u. 13--15,  Hungary}

\author{András Némethi}
\email{nemethi@renyi.hu}
\urladdr{\url{http://www.renyi.hu/~nemethi}}
\thanks{The second author is partially supported by NSF
      grant DMS-0304759, Marie Curie and OTKA grants.}
\address{Alfréd Rényi Institute of Mathematics\\
    1053 Budapest,  Reáltanoda u. 13--15,  Hungary}

\hypersetup{
  pdftitle={Invariants of Newton
    non-degenerate surface singularities},
  pdfauthor={Gábor Braun and András Némethi},
  pdfkeywords={hypersurface singularities, links of singularities,
    resolution graphs, Newton boundary, Newton polyhedrons},
  pdfsubject={2000 AMS Classification: Primary: 14J17, 14Q10;
  Secondary: 52B20}
}
\date{February 15, 2007}

\begin{document}
\begin{abstract}
  We recover the Newton diagram (modulo a natural ambiguity) from the link for
  any surface hypersurface singularity with non-degenerate Newton principal
  part whose link is a rational homology sphere. As a corollary, we show that
  the link determines the embedded topological type, the Milnor fibration, and
  the multiplicity of such a germ.  This proves (even a stronger version of)
  Zariski's Conjecture about the multiplicity for such a singularity.
\end{abstract}

\maketitle{}

\section{Introduction}
\label{sec:introduction}

In general, it is a rather challenging task to connect the analytic and
topological invariants of normal surface singularities. The program which aims
to recover different discrete analytic invariants from the abstract
topological type of the singularity (i.e.\ from the oriented homeomorphism
type of the link $K$, or from the resolution graph) can be considered as the
continuation of the work of Artin, Laufer, Tomari, S. S.-T.  Yau (and the
second author) about rational and elliptic singularities.  It includes the
efforts of Neumann and Wahl to recover the possible equations of the universal
abelian covers \cite{NW}, and the efforts of the second author and Nicolaescu
about the possible connections of the geometric genus with the Seiberg--Witten
invariants of the link \cite{NN}. See \cite{Nem} for a review of this program.

In order to have a chance for this program, one has to consider a topological
restriction (the weakest one for which we still hope for positive results
maybe that the link is a rational homology sphere), and a restriction about
the analytic type of the singularity, also.  By~\cite{SI}, the Gorenstein
condition is not sufficient.  We expect pathologies even for hypersurface
singularities.

For isolated hypersurface singularities a famous conjecture was formulated by
Zariski \cite{Za}, which predicts that the multiplicity is determined by the
\emph{embedded} topological type. For hypersurface germs with rational
homology sphere links, Mendris and the second author in \cite{MN} formulated
(and verified for suspension singularities) an even stronger conjecture,
namely that already the abstract link determines the embedded topological
type, the multiplicity and equivariant Hodge numbers (of the vanishing
cohomology).

The goal of the present article is to verify this stronger conjecture for
isolated singularities with non-degenerate Newton principal part. In fact, we
will prove that from the link (provided that it is a rational homology sphere)
one can recover the Newton boundary (up to a natural ambiguity, see
Theorem~\ref{th:1} below, and up to a permutation of coordinates), and hence
the equation of the germ (up to an equisingular deformation).  This is the
maximum what we can hope for.

The reader is invited to consult \cite{MR966191,Nem} for general facts about
singularities. §\ref{sec:sing-with-non} reviews the terminology and some
properties of germs with non-degenerate Newton principal part. In §\ref{2.6}
we define the equivalence relation $\sim$ of Newton boundaries characterizing the
above-mentioned ambiguities.  It may also be generated by the following
elementary step: two diagrams $\Gamma_1$ and $\Gamma_2$ are equivalent if both define
isolated singularities and $\Gamma_1\subset \Gamma_2$. (At the level of germs, this can be
described by a linear deformation.) Although the structure of an equivalence
class is not immediate from the definition, we define an easily recognizable
representative in every class, which we call the \emph{d-minimal}
representative.

In §\ref{OKAsect} we review Oka's algorithm which provides a possible
resolution graph $G(\Gamma)$ (or equivalently, a plumbing graph of the link) from
the Newton boundary $\Gamma$ \cite{MR894303}.  (Equivalent graphs provide plumbing
graphs related by blowing ups/downs, and hence determine the same link.) Our
main result says that Oka's algorithm can be essentially inverted:

\begin{theorem}\label{th:1} Assume that the Newton diagrams
  $\Gamma_1$ and $\Gamma_2$ determine isolated singularities with non-degenerate Newton
  principal part whose links are rational homology spheres.  Assume that the
  good minimal resolution graphs associated with $G(\Gamma_1)$ and $G(\Gamma_2)$ are
  isomorphic.  Then (up to a permutation of coordinates) $\Gamma_1\sim\Gamma_2$.  In
  particular, from the link $K$, one can identify the $\sim$-equivalence class of
  the Newton boundary (up to a permutation of coordinates) or, equivalently,
  the d-minimal representative of this class.
\end{theorem}

In fact, we prove an even stronger result: one can recover the corresponding
class of Newton diagrams (or its distinguished representative) already from
the \emph{orbifold diagram} $G^o$ associated with the good minimal resolution
graph. This diagram, a priori, contains less information then the resolution
graph, because it codifies only its shape and some subgraph-determinants,
see~\ref{sec:splice-diagram} for details.  (Although $G^o$ has a different
decoration, it is comparable with the `splice diagram' considered in
\cite{NW}.)

Since most of the invariants of the germs are stable under the deformations
defining the equivalence relation $\sim$ (see §\ref{stab}), one has the following
\begin{corollary}\label{cor:1}
  Let $f$ be an isolated germ with non-degenerate Newton principal part whose
  link is a rational homology sphere. Then the oriented topological type of
  its link determines completely its Milnor number, geometric genus, spectral
  numbers, multiplicity, and, finally, its embedded topological type.
\end{corollary}
Such a statement is highly non-trivial for any of the above invariants.  For
the history of the problem regarding the Milnor number and the geometric
genus, the reader is invited to consult \cite{Nem}. Here we emphasize only the
following:

\begin{itemize}
\item Regarding the embedded topological type, Corollary~\ref{cor:1} shows
  that if a rational homology sphere \(3\)-manifold can be embedded into $S^5$
  as the embedded link of an isolated hypersurface singularity with
  non-degenerate Newton principal part, then this embedding is \emph{unique}.
  (Notice the huge difference to the case of plane curves, and also to the
  higher dimensional case, where already the Brieskorn singularities provide a
  big variety of embeddings $S^{2n-1}\subset S^{2n+1}$, $n\neq2$.)

\item Such a link can be realized by a germ $f$ with non-degenerate Newton
  principal part in an essentially unique way, i.e.\ up to a sequence of
  linear $\mu$-constant deformations (corresponding to $\sim$) and permutation of
  coordinates, see~\ref{okaoka}(\ref{item:47}).
\end{itemize}

Regarding the main theorem, some more comments are in order.

\begin{itemize}
\item The assumption that the link is a rational homology sphere is necessary:
  the germs $\{z_1^a+z_2^b+z_3^c=0\}$ with exponents \((3,7,21)\) and
  \((4,5,20)\) share the same minimal resolution graph.
\item The proof of \ref{th:1} is, in fact, a \emph{constructive algorithm}
  which provides the d-minimal representatives of the corresponding class of
  diagrams from the orbifold diagram $G^o$.
\end{itemize}
Hence, one may check effectively whether an arbitrary resolution graph can be
realized by a hypersurface singularity with non-degenerate Newton principal
part.  Indeed, if one runs our algorithm and it fails, then it is definitely
not of this type. If the algorithm goes through and provides some candidate
for a Newton diagram, then one has to compute the graph (orbifold diagram) of
this candidate (by Oka's procedure) and compare with the initial one. If they
agree then the answer is yes; if they are different, the answer again is no
(this may happen since our algorithm uses only a part of the information of
$G^o$).

E.g., one can check that the following resolution graph cannot be realized by
an isolated singularity with non-degenerate Newton principal part (although it
can be realized by a suspension $\{z_3^2+g(z_1,z_2)=0\}$, where $g$ is an
irreducible plane curve singularity with Newton pairs \((2,3)\) and
\((1,3)\)).

\begin{figure}[h]
  \centering
  \begin{picture}(145,60)(-10,-10) \put(40,43){\makebox(0,0){$-7$}}
    \put(0,43){\makebox(0,0){$-3$}} \put(55,0){\makebox(0,0){$-3$}}
    \put(80,43){\makebox(0,0){$-1$}} \put(95,0){\makebox(0,0){$-3$}}
    \put(120,43){\makebox(0,0){$-2$}} \put(40,0){\circle*{5}}
    \put(40,33){\circle*{5}} \put(0,33){\circle*{5}} \put(80,33){\circle*{5}}
    \put(120,33){\circle*{5}} \put(80,0){\circle*{5}}
    \put(80,33){\line(0,-1){33}} \put(0,33){\line(1,0){120}}
    \put(40,33){\line(0,-1){33}}
  \end{picture}
  \caption{A resolution graph not coming from a singularity with
    non-degenerate Newton principal part}
  \label{fig:Newton-degenerate}
\end{figure}
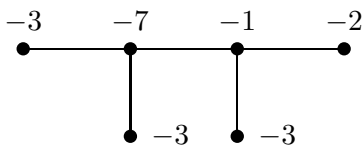
We mention that, in general, there is no procedure to decide whether a graph
is the resolution graph of a hypersurface isolated singularity (this is one of
the open problems asked by Laufer \cite[p.~122]{Open}; for suspension
singularities, it is solved in \cite{MN}).

\section{Singularities with non-degenerate Newton principal part}
\label{sec:sing-with-non}

\subsection{The Newton boundary\texorpdfstring{ \cite{MR0419433}}{}.
  Criterion for isolated singularities.}

\begin{para}\label{2.1} For any set  $S\subset\setN^3$ denote by $\Gamma_+(S)
  \subset \setR^3$ the convex closure of $\bigcup_{p\in S}(p+\setR_+^3)$. We call the \(1\)-faces of
  any polytope \emph{edges}, and \emph{face} will simply mean a \(2\)-face.
  The collection of all boundary faces of $\Gamma_+(S)$ is denoted by $\calf$.  The
  set of \emph{compact} faces of \(\Gamma_+(S)\) is denoted by $\calf_c$.  By
  definition, the \emph{Newton boundary} (or \emph{diagram}) \(\Gamma(S)\)
  associated with $S$ is the union of compact boundary faces of $\Gamma_+(S)$. Let
  $\partial\Gamma$ denote the union of those edges of $\Gamma(S)$ which are not intersection of
  two faces of $\Gamma(S)$.  Let $\Gamma_-(S)$ denote the cone with base \(\Gamma(S)\) and
  vertex \(0\).

  Let $f\colon (\setC^3,0) \to (\setC,0)$ be an analytic function germ defined by a
  convergent power series $\sum_p a_p z^p$ (where $p=(p_1,p_2,p_3)$ and
  $z^p=z_1^{p_1}z_2^{p_2}z_3^{p_3}$). By definition, the \emph{Newton
    boundary} $\Gamma(f)$ of $f$ is $\Gamma(\supp(f))$, where $\supp(f)$ is the support
  $\{p : a_p\neq0\}$ of $f$, and we write $\Gamma_-(f)$ for $\Gamma_-(\supp(f))$. The
  \emph{Newton principal part} of $f$ is $\sum_{p\in \Gamma(f)} a_p z^p$. Similarly, for
  any $q$-face $\bigtriangleup $ of $\Gamma(f)$ (of any dimension $q$), set $f_\bigtriangleup(z)\coloneqq
  \sum_{p\in\bigtriangleup}a_p z^p$. We say that $f$ is \emph{non-degenerate} on $\bigtriangleup$ if the
  system of equations $\partial f_\bigtriangleup/\partial z_1=\partial f_\bigtriangleup/\partial z_2 =\partial f_\bigtriangleup/\partial z_3=0$ has no solution
  in $(\setC^*)^3$. When $f$ is non-degenerate on every $q$-face of $\Gamma(f)$, we say
  (after Kouchnirenko \cite{MR0419433}) that $f$ has a \emph{non-degenerate
    Newton principal part}.  The diagram $\Gamma(f)$ and the function $f$ are
  called \emph{convenient} if $\Gamma(f)$ intersects all the coordinate
  axes.\end{para}
\begin{para}\label{eq:1}
  In this article we will assume that $f$ is singular, i.e.\ $\partial f(0)=0$.
\end{para}
\begin{para}\label{2.2} If we fix a Newton boundary $\Gamma$ (i.e.\
  $\Gamma=\Gamma(S)$ for some $S$), then the set of coefficients $\{a_p : p\in \Gamma\}$ for
  which $f(z)=\sum_{p\in \Gamma} a_p z^p$ is Newton non-degenerate (as its own principal
  part) form a non-empty Zariski open set (cf.\ \cite[1.10(iii)]{MR0419433}).
  Nevertheless, even for generic coefficients ${\{a_p\}}_{p\in\Gamma}$, the germ
  $f=\sum_{p\in \Gamma} a_p z^p$ (or any $f$ with $\Gamma(f)=\Gamma$), in general, does not define
  an isolated singularity.  The germ $f$ (with generic ${\{a_p\}}_{p\in \Gamma}$)
  defines an \emph{isolated} singularity if and only if $\Gamma$ satisfies the next
  additional properties (\cite[1.13(ii)]{MR0419433}):
  \begin{fakeequation}\label{eq:2}   
      \begin{itemize}
      \item \(\{(0,0,0),(0,0,1),(0,1,0),(1,0,0)\}\cap \Gamma=\emptyset\) (cf.~\eqref{eq:1}),
      \item the diagram $\Gamma$ has a vertex on every coordinate plane, and
      \item for every coordinate axis, $\Gamma$ has a vertex at most \(1\) far from
        the axis.
      \end{itemize}
  \end{fakeequation}
  E.g., a convenient $f$ with generic coefficients defines an isolated
  singularity.
\end{para}

\begin{example}\label{ex:100} Notice that \eqref{eq:2}
  cannot be satisfied by one vertex. Moreover, if $\Gamma$ satisfies \eqref{eq:2}
  and has no faces then (modulo a permutation of the coordinates) it is the
  segment \([(0,1,1),(n,0,0)]\) for some $n\geq 2$.
\end{example}
\begin{remark}\label{rem:uj} Assume that $\Gamma$ is not an
  edge. Then \eqref{eq:2} implies that every edge of $\partial \Gamma$ should lie either
  on a coordinate plane or be (after permuting coordinates) of the form
  \(AB=[(a,0,c),(0,1,b)]\) with $a>0$ and $b+c>0$. The number of edges of
  second type coincides with the number of coordinate axes not intersected by
  $\Gamma$. (Indeed, assume that the \(z\sb{3}\) axis does not meet $\Gamma$. Project
  $\Gamma$ to the \(z\sb{1} z\sb{2}\) plane by $\psi(z_1,z_2,z_3)=(z_1,z_2)$. Then, by
  \eqref{eq:2}, the boundary of $\psi(\Gamma)$ contains an edge of type
  $[(a,0),(0,1)]$.) \end{remark}

\begin{para} If one tries to analyze the invariants of a germ in terms of its
  Newton diagram (see e.g.\ the references cited in §\ref{2.3}), one
  inevitably faces the arithmetical properties of integral polytopes. In
  Appendix~\ref{arith}, we collect those which will be used in the body of the
  paper. The relevant notations and terminologies are listed below:
\end{para}
\begin{para}[Notations/Definitions]\label{not} Fix a Newton
  diagram. Set $\bigtriangleup\in\calf$. Let $\bigtriangledown\in\calf$ be an adjacent face with a common
  (compact) edge $AB\coloneqq \bigtriangleup\cap \bigtriangledown$. Then one defines:

  \begin{center}
    \begin{tabularx}{0.9\linewidth}{>{\(}c<{\)}X}
      \normalvector{\bigtriangleup}& the \emph{normal vector} of \(\bigtriangleup\), i.e.\ the primitive
      integral
      vector with non-negative entries, normal to \(\bigtriangleup\),\\
      \multiplicity{\bigtriangleup}{\bigtriangledown}&
      the number of components of $AB\setminus \setN^3_{>0}$,\\
      \innerdeterminant{\bigtriangleup}{\bigtriangledown}& the \emph{determinant} of \(\normalvector{\bigtriangleup}\)
      and \(\normalvector{\bigtriangledown}\), namely, the greatest common divisor of the
      entries of the cross product \(\normalvector{\bigtriangleup} × \normalvector{\bigtriangledown}\),
      (\(\innerdeterminant{\bigtriangleup}{\bigtriangledown}\geq 1\)),\\
      \coordinatevector{1}, \coordinatevector{2}, \coordinatevector{3}
      &the three coordinate normal vectors.\\
    \end{tabularx}
  \end{center}
  The number $ \innerdeterminant{\bigtriangleup}{\bigtriangledown}$ is also called the determinant of the
  edge $AB$.  Since it depends only on the corresponding normal vectors,
  sometimes we put the normal vectors in the index instead of the faces. E.g.,
  if $\normalvector{\bigtriangledown}=\coordinatevector{i}$ and \(\bigtriangleup\in \calf_c\), then we may
  also write $\innerdeterminant{\bigtriangleup}{\coordinatevector{i}}$ for
  $\innerdeterminant{\bigtriangleup}{\bigtriangledown}$.  The number
  $\multiplicity{\bigtriangleup}{\coordinatevector{i}}$ has a similar meaning. In fact,
  with the notation $\normalvector{\bigtriangleup}=(a_1,a_2,a_3)$, one has:
  \begin{equation}
    \label{eq:ijk}
    \innerdeterminant{\bigtriangleup}{\coordinatevector{i}}=
    \mygcd{a_{j},
      a_{k}},\quad
    \text{where $\{i,j,k\}=\{1,2,3\}$}.
  \end{equation}
  Similarly, for any lattice polygon $\bigtriangleup$, the vector \(\normalvector{\bigtriangleup}\)
  denotes the primitive integral vector normal to \(\bigtriangleup\) (well-defined up to a
  sign).  The \emph{combinatorial area}, by definition
  (cf.~\cite[(6.2)]{MR894303}), is
  \begin{equation}
    \label{eq:3}
    \combinatorialarea{\bigtriangleup} \coloneqq 2 \#\{\text{inner lattice points}\}
    + \#\{\text{border lattice  points}\} - 2.
  \end{equation}
  Clearly, \(\combinatorialarea{\bigtriangleup}\) is additive.  The face \(\bigtriangleup\) is called
  \emph{empty} if its only lattice points are its vertices.
\end{para}

\subsection{Some discrete invariants determined from the Newton
  boundary.}\label{2.3}

If $f$ defines an isolated singularity and has a non-degenerate Newton
principal part, then its Newton boundary $\Gamma(f)$ determines almost all its
discrete analytic and embedded topological invariants. E.g.:
\begin{enumerate}[(a)]
\item the \emph{Milnor number} $\mu(f)$ of $f$ is given by Kouchnirenko
  \cite{MR0419433}.  For any $\Gamma$ let $V_3$ be the \(3\)-dimensional volume of
  $\Gamma_-$, and for $1\leq q\leq 2$, let $V_q$ be the sum of the $q$-dimensional
  volumes of all the intersections of $\Gamma_-$ with $q$-dimensional coordinate
  planes. Set $\nu(\Gamma)\coloneqq 6V_3-2V_2+V_1-1$. Then, by \cite{MR0419433}, the
  Milnor number $\mu(f)$ of any \emph{convenient} germ $f$ with non-degenerate
  Newton principal part is given combinatorially via $\Gamma(f)$ by:
  \begin{equation}
    \label{eq:4}
    \mu(f) =\nu(\Gamma(f)).
  \end{equation}
  In fact, the same formula is valid for non-convenient isolated singularities
  as well. Indeed, assume e.g.\ that the diagram $\Gamma(f)$ does not intersect the
  \(z\sb{3}\) axis, and let \(AB\) be an edge as in \ref{rem:uj}. Then the
  deformation $f_d\coloneqq f+tz_3^d$ with $d\geq \mu(f)+2$ has a uniform stable
  radius for the Milnor fibration \cite{MR535090}, hence $\mu(f)=\mu(f_d)$.
  Moreover, $\Gamma_-(f_d)=\Gamma_-(f)\cup W_d$, where $W_d$ is the \(3\)-simplex with
  vertices \(0\), \(A\), \(B\) and \((0,0,d)\). Since
  $(6V_3-2V_2+V_1)(W_d)=0$, one gets that $\nu(\Gamma(f_d))=\nu(\Gamma(f))$.

  (Since $f$ is finitely determined, $f$ and $f_d$ are right-equivalent for
  $d\gg 0$ and their other invariants listed in this subsection agree too.
  Hence, being convenient, in many cases present in the literature, is not
  really essential for us, see also \cite{Wall}.)

\item the \emph{characteristic polynomial} of the algebraic monodromy is
  determined in \cite{MR0424806}; the \emph{geometric genus} of the surface
  singularity $(\{f=0\},0)$ is given by $\#(\Gamma_-(f)\cap \setN_{>0}^3)$,
  cf.~\cite{0461.14009,MR954149}; the set of \emph{spectral numbers} (or
  characteristic exponents) is computed in \cite{Da,MR954149,ST,VK}; the
  \emph{multiplicity} of $f$ by $\min_{p\in \Gamma(f)} \sum p_i$;
\item the \emph{embedded topological type} and the \emph{Milnor fibration} of
  $f$ (with its homological `package' including the Seifert form) is
  determined from $\Gamma(f)$ uniquely by \cite[2.1]{MR535090};
\item an explicit construction of the \emph{dual resolution graph} $G(f)$ of
  the surface singularity $(\{f=0\},0)$ is given in \cite{MR894303} (we review
  this in §\ref{OKAsect}).
\end{enumerate}

\subsection{The structure of Newton polytopes in the case of rational homology
  sphere links}\label{sec:poss-newt-diagr}

\begin{para}\label{rrhs}
  An important assumption of the main result of the present article is that
  the link $K(f)$ of $f$ is a \emph{rational homology sphere}, i.e.\ $
  H_1(K(f),\setQ)=0$. This additional assumption (besides \eqref{eq:2}, which says
  that $f$ with non-degenerate Newton principal part is an isolated
  singularity) imposes serious restrictions on the Newton boundary $\Gamma(f)$,
  cf.~\cite{MR954149}:
  \begin{equation}\label{eq:rhs}
    \text{$K(f)$ is a rational homology sphere}  \iff
    \Gamma(f)\cap \setN_{>0}^3=\emptyset.
  \end{equation}
\end{para}
In this subsection \emph{we assume that $\Gamma(f)$ satisfies these two
  restrictions, namely \eqref{eq:2} and \eqref{eq:rhs}}.  Our goal is to
derive the structure theorem~\ref{prop:1} for Newton diagrams.

We fix a diagram $\Gamma$. We start by classifying the non-triangular faces:

\begin{lemma} \label{lem:1} If a face of $\Gamma$ is not a triangle then it is a
  trapezoid.  By permuting coordinates, its vertices are: \(A=(p,0,n)\),
  \(B=(0,q,n)\), \(C=(r_1,r_2+tq,0)\) and \(D=(r_1+tp,r_2,0)\), where \(p,q >
  0\), \(\mygcd{p,q}=1\), \(t \geq 1\) and \(r_1,r_2 \geq 0\).  The only side which
  can have inner lattice points is the base lying on the \(z\sb{1} z\sb{2}\)
  plane (with \(t-1\) of them).
\end{lemma}

\begin{proof} The idea of the proof is the following: if a lattice polygon
\(\bigtriangleup\) is not a triangle or a trapezoid, then there exists a parallelogram in
\(\bigtriangleup\) with three vertices on the boundary of \(\bigtriangleup\) and one in its interior,
which contradicts \eqref{eq:rhs}. The details are left to the reader.
\end{proof}

\begin{terminology}\label{term1}
  The \emph{edges of a trapezoid} have asymmetric roles.  For future reference
  we give names to them. The \emph{bottom} edge always lies on a coordinate
  plane. If two (ore more) edges lie on coordinate planes, the bottom edge is
  the one which has internal lattice points, if such exists.  Otherwise, we
  choose one of them arbitrarily.

  Opposite to the bottom edge lies the \emph{top} edge, and the others are
  called \emph{side} edges.
\end{terminology}
\begin{termdis}\label{term2}
  An edge \emph{crosses}, say, the \(z\sb{3}\) axis if it is of the form
  \([(p,0,a),(0,q,b)]\), where $p>0$, $q>0$, and \(a+b>0\). There are two
  types of edges on $\Gamma$: those lying on a coordinate plane and those crossing
  a coordinate axis.

  While edges of the first type do not `cut' $\Gamma$, edges of the second type
  usually cut $\Gamma$ into two non-empty parts, one of which has a particularly
  simple structure. In order to see this, project $\setR_{\geq 0}^3\setminus 0$ from the
  origin to the triangle $\setT\coloneqq \{z_1+z_2+z_3=1 : z_i\geq 0 \
  (i=1,2,3)\}$. The restriction $\phi\colon \Gamma\to \setT$ is one-to-one and preserves
  segments.  An edge lying on a coordinate plane projects into $\partial \setT$,
  while a crossing edge projects into a segment with only its end points on $\partial
  \setT$ and cutting $\setT$ into two parts such that at least one of them,
  say $\setT_{0}$, is a triangle.  By~\ref{lem:1}, the projection of a
  trapezoid hits the interior of all the sides of $\setT$, hence
  $\phi^{-1}(\setT_{0})$ may contain only a `sequence of triangles'. Therefore,
  one has:
\end{termdis}
\begin{lemma}\label{lem:2}
  An edge of $\Gamma$ crossing (say) the \(z\sb{3}\) axis, which is not on $\partial\Gamma$,
  cuts $\Gamma$ into two non-empty parts.  Consider the plane $\pi$ formed by the
  edge and the origin. Then that part of $\Gamma$, which is on the same side of $\pi$
  as the positive \(z\sb{3}\) axis, consists only of triangular faces with
  vertices lying on the $z_1z_3$ and $z_2z_3$ planes.  They form a sequence
  $\bigtriangleup_1,\dotsc, \bigtriangleup_k$; where $\bigtriangleup_i$ is adjacent with $\bigtriangleup_{i+1}$, (and these are
  the only adjacent relations).  \end{lemma}

\begin{cordef}\label{cor:2}
  Fix a coordinate axis.

  First, assume that there is at least one triangular face whose vertices are
  on the two coordinate planes adjacent to the axis. Then the collection of
  such triangular faces form a sequence as in \ref{lem:2}, and their union is
  called the \emph{arm} of the diagram in the direction of that axis.  The arm
  also contains all the crossing edges whose vertices lie on the two
  coordinate planes.  Let the \emph{hand} be the triangle of the arm which is
  nearest to the axis (in the $\phi$-projection, say).  Let the \emph{shoulder}
  be the crossing edge of the arm which is most distant from the axis (in the
  same sense).

  Next, assume that there is no triangular face whose vertices are on these
  two coordinate planes.  Then we distinguish two cases:
  \begin{enumerate}[(a)]
  \item\label{item:1} If there exists a crossing edge of the coordinate axis,
    then it is unique; in this case we say that the arm in that direction is
    \emph{degenerate}, and the \emph{degenerate arm} (and its shoulder too) is
    this unique crossing edge.
  \item\label{item:2} If there is no crossing edge either, then we say that
    there is no arm in the direction of the axis.
  \end{enumerate}
\end{cordef}

\begin{terminology}\label{term3}
  A triangular face of $\Gamma$ is called \emph{central} if its vertices are not
  situated on the union of two coordinate planes. A face of $\Gamma$ is called
  \emph{central} if it either is a central triangle or it is a trapezoid. An
  edge of $\Gamma$ is \emph{central} if (modulo a permutation of the coordinates)
  it has the form \([(0,0,a),(p,q,0)]\).

  Using the projection $\phi\colon \Gamma\to\setT$, one may easily verify:

\end{terminology}

\begin{lemma}\label{lem:cf}
  $\Gamma$ has at most one central face. $\Gamma$ has a central face if and only if it
  has no central edge.
\end{lemma}

These facts can be summarized in the next result on \emph{structure of Newton
  diagrams}:

\begin{proposition}\label{prop:1}
  Every Newton diagram $\Gamma$ (which satisfies \eqref{eq:2} and \eqref{eq:rhs})
  sits in exactly one of the three disjoint families characterized as follows:
  \begin{enumerate}[(1)]
  \item\label{item:3} $\Gamma$ has a unique central trapezoid with at most \(3\)
    disjoint (possibly degenerate) arms. The arms correspond to those sides of
    the trapezoid which are crossing edges.

  \item\label{item:4} $\Gamma$ has a unique central triangle with \(3\) disjoint
    (possibly degenerate) arms.

  \item\label{item:5} $\Gamma$ has (at least one) central edge.
  \end{enumerate}

  Moreover, if \(\Gamma\) has a central edge, then there are two cases.  If $\Gamma$ has
  only one face, this face is triangular with all vertices on coordinate axes,
  then all edges are central.  Otherwise, all central edges have a common
  intersection point (say $P$) sitting on a coordinate axis; and the diagram
  has two (possibly degenerate) arms in the direction of the other two axes.
  The arms may overlap each other, i.e.\ have common triangles.  $P$ is a
  vertex of all the triangles in the intersection of the arms, and all those
  edges of these triangles which contain $P$ are central (and these are all
  the central edges).
\end{proposition}

\section{Equivalent Newton boundaries. Deformations.}\label{2.6}

\subsection{The equivalence relation}\label{eqrel}
\begin{para}
  Our aim is to recover the Newton boundary (up to a permutation of
  coordinates) of an isolated singularity with non-degenerate Newton principal
  part from the link $K(f)$, provided that $K(f)$ is a rational homology
  sphere. Strictly speaking, this is not possible: one can easily construct
  pairs of such germs having identical links but different boundaries. E.g.,
  take an isolated non-convenient germ $f$ and $f_d= f+\sum z_i^d$ with $d\gg 0$.
  This motivates to define a natural equivalence relation of Newton
  boundaries. By definition, it will be generated by two combinatorial
  `steps'.
\end{para}
\begin{para}\label{moves}
  Fix a Newton boundary $\Gamma=\Gamma(S)$ which satisfies~\eqref{eq:2}. Let $AB$ be an
  edge of $\partial\Gamma$ which is not contained in any coordinate plane.
  By~\ref{rem:uj}, up to a permutation of coordinates, $A=(a,0,c)$ (with
  $a>0$) and $B=(0,1,b)$.
\end{para}
\begin{enumerate}[Move 1.]
\item\label{item:6} We add a new vertex $C=(a',0,c')$ to $\Gamma$ in such a way
  that $\Gamma(S\cup C) =\Gamma(S)\cup \trianglebyvertices{A}{B}{C}$.  Here
  $\trianglebyvertices{A}{B}{C}$, the \(2\)-simplex spanned by the points
  \(A\), \(B\), \(C\) appears as a new face.  (In particular, $0\leq a'<a$ and
  $c'$ must be sufficiently large.)
\item\label{item:7} Assume that $AB$ is in the face $\bigtriangleup$ whose supporting plane
  is $H$. The line through $AB$ cuts out the open semi-plane $H_+$ of $H$
  which does not contain $\bigtriangleup$.  Set $S'\coloneqq H_+\cap \setN^3$.  Then by adding a
  non-empty subset $S''$ of $S'$ to $S$, we create a new Newton boundary $\Gamma(S\cup
  S'')$.  By this move, all faces of $\Gamma(S)$ are unmodified, except $\bigtriangleup$, which
  is replaced by a larger face containing $\bigtriangleup$. \end{enumerate}
\begin{definition}\label{2.7}
  We denote Move~\ref{item:6} and Move~\ref{item:7} by $\ep$ and $\kp$,
  respectively.  We denote their inverses by $\en$ and $\kn$, respectively.
  The segment $AB$ will be called the \emph{axis} of the corresponding move.

  Two Newton diagrams $\Gamma_1$ and $\Gamma_2$, both satisfying~\eqref{eq:2}, are
  \emph{equivalent} (and we write $\Gamma_1\sim \Gamma_2$), if they can be connected by a
  sequence of elementary moves ($\epn$ or $\kpn$), such that all the
  intermediate Newton boundaries satisfy \eqref{eq:2} as well.
\end{definition}

\begin{example}\label{elso}
  Using \ref{rem:uj} and induction, one can show that if $\Gamma_1$ and $\Gamma_2$ are
  Newton diagrams, $\Gamma_1\subset \Gamma_2$, both satisfying~\eqref{eq:2}, then they are
  equivalent. In fact, the inclusion of Newton boundaries with~\eqref{eq:2}
  generates the same equivalence relation.
\end{example}

\begin{example}\label{ex:101}
  The segments \([(0,1,1),(n,0,0)]\) and \([(1,0,1),(0,n,0)]\) (considered as
  diagrams) are equivalent. Indeed, add to \(\Gamma_1= [(0,1,1),(n,0,0)]\) the
  vertex $(1,0,1)$ (by $\ep$), then add $(0,n,0)$ (by $\kp$), then remove the
  end points of $\Gamma_1$ (cf.~\ref{elso}).
\end{example}

\begin{para}\label{linear}
  Sometimes it is more convenient to specify the deformation of the
  corresponding germs instead of the modification of Newton diagrams: adding a
  new vertex $p$ to $S$ translates into adding a new monomial $ta_pz^p$ to
  $f$, with $t\in[0,\epsilon]$ a deformation parameter. (The fact that these
  deformations are \emph{linear} in $t$ is crucial in the proof
  of~\ref{2.8}(\ref{item:10})).
\end{para}
\begin{example}\label{ex:moving-triangle}
  The number of `essential' deformation parameters can be as large as we wish.
  E.g., for $m,n\gg 0$, all the different Newton diagrams associated with the
  family
  \begin{equation*}
    z_3(z_1^p+z_2^q+z_3^r)+ \sum_i t_iz_1^{m-ip}z_2^{n+iq}
    \quad (m-ip\geq 0, n+iq\geq0)
  \end{equation*}
  satisfy \eqref{eq:2}, and are equivalent (via repeated $\kpn$) as soon as
  $\sum_i \lvert t_i \rvert > 0$.  We call the `ambiguity' of the choice of the
  monomials $z_1^{m-ip}z_2^{n+iq}$ the \emph{moving triangle ambiguity}.

  More generally, a \emph{moving triangle} of a Newton diagram $\Gamma$ is a
  triangular face with vertices: \(P\coloneqq(p,0,1)\), \(Q\coloneqq(0,q,1)\)
  and \(R \coloneqq (m,n,0)\), where the edge $PQ$ is in some other face as
  well.  Consider the line through \(R\) parallel to $PQ$.  Then (the
  \emph{moving vertex}) \(R\) can be replaced by any of the lattice points
  \(S\) on this line with non-negative coordinates (or any collection of
  them).  If $\Gamma$ satisfies \eqref{eq:rhs}, then $\mygcd{p,q}=1$, and by
  \eqref{eq:19} $ \normalvector{\bigtriangleup}= (a_1,a_2,a_3)=(q,p,mq+np-pq)$.  Therefore,
  one has:
  \begin{equation}\label{eq:movcor}
    \begin{alignedat}{4} p &\mid m  &&\iff
      &a_{2} &\mid a_{3}
      &&\iff \text{$R$ can be replaced by a point on the \(z\sb{2}\) axis},\\
      q &\mid n &&\iff
      &a_{1} &\mid a_{3}
      &&\iff \text{$R$ can be replaced by a point on the $z_1$
        axis}.
    \end{alignedat}\end{equation}
\end{example}

\subsection{Stability of the invariants under the deformations.}\label{stab}

\begin{proposition}\label{2.8} Consider two isolated singularities with non-degenerate Newton principal parts whose Newton boundaries are equivalent in the
  sense of \ref{2.7}. Then the following invariants associated with these
  germs are the same:
  \begin{enumerate}[(a)]
  \item\label{item:8} the Milnor number $\mu$;
  \item\label{item:9} the link $K$;
  \item\label{item:10} more generally, the embedded topological type;
  \item\label{item:11} the spectral numbers (in particular, the geometric
    genus); the equivariant Hodge numbers;
  \item\label{item:12} the multiplicity.
  \end{enumerate}
  Moreover, a deformation associated with $\ep$ or $\kp$ admits a weak
  simultaneous resolution.
\end{proposition}

\begin{proof} First of all, (\ref{item:8}) can be easily verified by direct
computation (left to the reader) by Kouchnirenko's formula \eqref{eq:4}.
Item~(\ref{item:9}) can also be checked directly from Oka's algorithm
\cite{MR894303} (§\ref{OKAsect} here), and (\ref{item:12}) is also elementary.
But there are also (more) conceptual short-cuts: The existence of a weak
simultaneous resolution follows from a result of Oka \cite{MR1000603} (after
we add some high degree monomials in the non-convenient case, and we notice
that our moves are `negligible truncations' in the sense of Oka), which
implies (\ref{item:9}) by a result of Laufer \cite{LauferW}.  For
(\ref{item:11}) one can use Varchenko's result \cite{Var}, which says that the
spectrum is constant under a $\mu$-constant deformation. Notice also that the
geometric genus is the number of spectral numbers in the interval $(0,1]$.
Finally, a $\mu$-constant $(f+tg)$-type deformation (cf.~\ref{linear}) is
topological trivial by a result of Parusiński \cite{Par} (proving
(\ref{item:10})), and is equimultiple (e.g.) by Trotman \cite{Tr}.
\end{proof}

\begin{remark} By similar proof as in \cite{Stsc} (valid for the spectrum),
  one can show that \emph{the set of spectral pairs} (equivalently, the
  equivariant Hodge numbers) of $f$ are also determined by $\Gamma$, and are stable
  with respect to the $\sim$-deformation.  Cf.\ also with \cite{Da}.
\end{remark}

\begin{corollary}\label{cor:3}
  Fix a Newton diagram $\Gamma$ which satisfies \eqref{eq:2}. Then the following
  facts are equivalent:
  \begin{enumerate}[(a)]
  \item\label{item:13} $\{(0,1,1),(1,0,1),(1,1,0)\}\cap \Gamma\neq\emptyset;$

  \item\label{item:14} $\Gamma$ is equivalent to a diagram which has no \(2\)
    dimensional faces;

  \item\label{item:15} $\Gamma$ is equivalent to the segment-diagram
    \([(0,1,1),(n,0,0)]\) for some $n\geq 2$;

  \item\label{item:16} $f_\Gamma(z)\coloneqq \sum_{p\in\Gamma} a_pz^p$ (with generic
    coefficients ${\{a_p\}}_p$) is an $A_{n-1}$ singularity (the unique
    hypersurface cyclic quotient singularity with $\mu=n-1$) for some $n\geq 2$.

  \item\label{item:17} The minimal dual resolution graph of $\{f_\Gamma=0\}$ is a
    string (with determinant $n$).
  \end{enumerate}
  In fact, the integers $n$ in (\ref{item:15}), (\ref{item:16}) and
  (\ref{item:17}) are equal.
\end{corollary}

\begin{proof}
(\ref{item:14})$\implies$(\ref{item:15}) follows from \ref{ex:100} and
\ref{ex:101}.  The implication (\ref{item:15})$\implies$(\ref{item:14}) is
clear. For (\ref{item:13})$\implies$(\ref{item:15}), using \ref{ex:101}, it is
enough to prove that if $(0,1,1)\in \Gamma$ then \(\Gamma\sim [(0,1,1),(n,0,0)]\) for some
$n$. If $\Gamma$ intersects the \(z\sb{3}\) axis at some point $(n,0,0)$, then
$[(0,1,1),(n,0,0)]\subset\Gamma$ and one may use \ref{elso}. Otherwise, one considers,
like in \ref{rem:uj}, the projection $(z_1,z_2,z_3)\mapsto (z_1,z_2) $ restricted to
$\Gamma$. By \eqref{eq:2} there is at least one edge whose projection has the form
$[(a,0),(0,1)]$ (up to a permutation). Consider the edge-projection of this
type which is closest to $(1,1)$, let its preimage in $\Gamma$ be \([
(a,0,c),(0,1,b)]\). This choice guarantees that $\Gamma$ contains the triangular
face \(\bigtriangleup\) with vertices \((a,0,c)\), \((0,1,b)\), \((1,1,0)\).  Since \(\bigtriangleup\),
as a diagram, satisfies \eqref{eq:2}, $\Gamma\sim\bigtriangleup$ by \ref{elso}. By $\kp$ one can
add to \(\bigtriangleup\) the vertex $(ab+c,0,0)$, and apply again \ref{elso} to show that
\(\bigtriangleup\sim [(0,1,1),(n,0,0)]\) with $n=ab+c$.

Next, notice that $a_1z_1^n+a_2z_2z_3$ ($a_1,a_2\neq0$) defines an $A_{n-1}$
singularity. Hence (\ref{item:15})$\implies$(\ref{item:16}) follows from
\ref{2.8}(\ref{item:9}), since the $A_{n-1}$ singularity is characterized by
the fact that its link is the lens space $L(n,n-1)$. For
(\ref{item:16})$\implies$(\ref{item:13}) one uses that the quadratic part of
the Taylor-expansion of an $A_{n-1}$ singularity $f_\Gamma$ (in any coordinate
system) has rank at least two.

(\ref{item:16})$\iff$(\ref{item:17}) follows from the fact that the $A_{n-1}$
singularities are the only hypersurface singularities whose minimal resolution
graphs are strings.
\end{proof}

\begin{para}
  Not all discrete analytic invariants of the germs remain constant under the
  above equivalence relation. The following example was provided by J. F. de
  Borbadilla, A. Melle-Hernández and I. Luengo (private communication), in
  which Teissier's invariant $\mu^*$ jumps.
\end{para}
\begin{example}
  Consider the deformation $f_t=z_3^3+z_2^4z_1+z_1^{10}+tz_2^3z_3$, which
  corresponds to Move~\ref{item:6}, hence $f_1$ and $f_0$ are equivalent in
  the sense of \ref{2.7}.  But the Milnor numbers of the generic hyperplane
  sections are not the same: $\mu^{(2)}(f_1)=7$, while $\mu^{(2)}(f_0)=8$. In
  particular, by \cite{Laufer}, this deformation does not admit a strong
  simultaneous resolution.  Similar example was constructed by Briançon and
  Speder \cite{BS} (cf.\ also with \cite{MR1000603}); the main difference is
  that in the present case the stable link $K(f_1)=K(f_0)$ is a rational
  homology sphere. Notice also that $f_0$ is weighted homogeneous and
  $\deg(z_2^3z_3)>\deg(f_0)$.  (This example also shows that
  \cite[Question~13.12]{Nem} has a negative answer: i.e.\ for a deformation
  which admits a weak simultaneous resolution the existence of a strong
  simultaneous resolution is not guaranteed, even if the stable link is a
  rational homology sphere.)
\end{example}
\begin{example}
  The recent manuscript \cite[§4]{Artal2007} provides a \(\mu\)-constant
  deformation of singularities with non-degenerate Newton principal part and
  \(b_1(K)>0\) such that the homeomorphism of the tangent cone jumps,
  providing a counterexample to \cite[Conjecture B.]{Za}.  In fact, a
  counterexample also exists among rational homology spheres, e.g.\ the
  deformation \(f_t = z_1^3 z_2 + z_3^5 + z_2^{11} + t z_1^2 z_3^2\) of type
  Move~\ref{item:6}.  Then the homeomorphism type of \(\{ z_1^3 z_2 + t z_1^2
  z_3^2 = 0 \}\) is not constant.
\end{example}

\subsection{Distinguished representatives}\label{representatives}

\begin{para}\label{princ}
  In this subsection we assume that \emph{all our Newton diagrams satisfy
    \eqref{eq:2} and \eqref{eq:rhs}}.  It is preferable to have in each
  $\sim$-equivalence class a well-characterized and easily recognizable
  representative to work with. In its choice we are guided by the following
  principles (motivated by \ref{mingood}, which says that such a `minimal'
  diagram reflects better the \emph{minimal} resolution graph $G_{\min}(f)$ of
  the germ $f$):
  \begin{enumerate}[(a)]
  \item\label{item:18} the representative should have a minimal number of
    faces;
  \item\label{item:19} all the faces which cannot be eliminated by $\en$
    should be `minimized as much as possible' by $\kn $;
  \item\label{item:20} a representative may contain a trapezoid only if the
    trapezoid cannot be replaced by a triangle in its class.
  \end{enumerate}
  This motivates the following:
\end{para}
\begin{definition}\label{def:1}
  A diagram $\Gamma$ is called \emph{M1-minimal} if by the direct application of a
  move of type $\en$ one cannot eliminate any of its faces.
\end{definition}

Notice that at least one M1-minimal representative exists in any equivalence
class.

Also, one can decide the M1-minimality of a diagram by analyzing the lattice
points sitting on it, without any information about the other diagrams in its
class. But, exactly for this reason, the above definition does not exclude the
possibility that an M1-minimal diagram may have another diagram in its class
with less faces. In fact, this may occur:

\begin{example}\label{ex:105} (Cf.\ the proof of Corollary~\ref{cor:3}.)
  The diagram consisting of the unique triangular face with vertices
  $(a,0,c),(0,1,b),(1,1,0)$ is M1-minimal, but it is equivalent to the segment
  \([(1,1,0),(0,0,ab+c)]\).
\end{example}
The next lemma guarantees that this is the only pathological case when such a
phenomenon may occur. Below $\#(\Gamma)$ denotes the number of faces of $\Gamma$.

\begin{lemma}\label{lem:3}
  Fix a diagram $\Gamma$, which is not of type characterized by
  Corollary~\ref{cor:3}.
  \begin{enumerate}[(a)]
  \item\label{item:21} Then $\Gamma$ is M1-minimal if and only if for any $\Gamma'\sim \Gamma$
    one has $\#(\Gamma')\geq \#(\Gamma)$.
  \item\label{item:22} If $\Gamma$ and $\Gamma'$ are both M1-minimal and $\Gamma\sim\Gamma'$, then
    they can be connected by a sequence of diagrams related to each other only
    by moves $\kpn$. In particular, the set of supporting planes of the faces
    of the two diagrams are the same.
  \end{enumerate}
\end{lemma}
Notice that the assumption is essential for part (\ref{item:22}), too: see
e.g.\ the segments of~\ref{ex:101}.
\begin{proof}
A sequence of diagrams $\Gamma_1,\Gamma_2,\dotsc,\Gamma_k$ \emph{connects} $\Gamma_1$ and $\Gamma_k$ if
$\Gamma_i$ and $\Gamma_{i+1}$ (for $1\leq i\leq k-1$) are related by one of the moves $\epn$
or $\kpn$, denoted by $\Gamma_i\move{\jpn}\Gamma_{i+1}$. Our goal is to replace a given
sequence of diagrams connecting $\Gamma$ and $\Gamma'$ by another one which has the
additional property that all \(\en\) moves appear first.  For this, first we
analyze how one can modify two consecutive moves in a sequence, where the
second one is $\en$:

\begin{fact}
  $\Gamma_i\move{\jpn}\Gamma_{i+1}\move{\en}\Gamma_{i+2}$ can be replaced either by moves
  $\Gamma_i\move{\en}\Gamma_{i+1}'\move{\jpn}\Gamma_{i+2}$, or by a single move of type
  $\Gamma_i\move{\en}\Gamma_{i+2}$, or both moves can be eliminated, i.e.\
  $\Gamma_i=\Gamma_{i+2}$.
\end{fact}

Indeed, if the two moves operate on different faces of the diagram then they
can be performed in the reverse order with the same effect.  So we can assume
that the two moves operate on the same face \(\bigtriangleup\).

If, additionally, the two moves have the same axis $AB$ (see~\ref{2.7}), then
the moves eliminate one triangle from both sides of $AB$, hence $\Gamma_i\sim AB$
contradicting our assumption. Hence the two axes are different.  This can
occur only if $\bigtriangleup$ is a triangle and the composition of the two moves is the
removal of \(\bigtriangleup\), i.e.\ a move of type \(\en\).  This finishes the proof of
the fact.

The easy consequence of the fact is that if \(\Gamma\sim\Gamma'\) then they can be
connected by a sequence \(\Gamma_i\), in which all the \(\en\) moves appear first
(preceding the moves of other type).  In particular, if \(\Gamma\) is M1-minimal
then this sequence does not contain any \(\en\) moves, so the number of faces
\(\#(\Gamma_i)\) is non-decreasing along the sequence.  This proves the non-trivial
part of \ref{lem:3}(\ref{item:21}).  If \(\Gamma'\) is also M1-minimal then
(applying the above also for the reverse sequence) \(\#(\Gamma_i)\) must be
constant, i.e.\ the sequence does not contain any move of type \(\epn\),
finishing the proof of~(\ref{item:22}).
\end{proof}

\begin{definition}[Canonical and minimal representatives]\label{def:2}
  Fix the equivalence class of a diagram which does not satisfy \ref{cor:3},
  and consider all M1-minimal representatives. By \ref{lem:3}(\ref{item:22})
  they are related to each other by moves $\kpn$.  Clearly, this set has a
  \emph{unique maximal element with respect to $\kpn$} (or equivalently, with
  respect to the inclusion). This diagram will be called the \emph{canonical
    representative} of the class. It can be easily recognized: it is
  M1-minimal, and all its faces are as large as possible.

  The canonical representative satisfies the principle (\ref{item:18}) of
  \ref{princ}, but not (\ref{item:19}). For (\ref{item:19}), we would need the
  unique \emph{minimal element with respect to $\kpn$} of all M1-minimal
  representatives; but such an element, in general, does not exist.
  Nevertheless, we consider the \emph{set of minimal elements} (diagrams which
  cannot be reduced by $\kn$) of all M1-minimal representatives. We call these
  representatives \emph{minimal}. By \ref{lem:3}, these are those
  representatives which cannot be reduced by any move $\jn$.
\end{definition}

\begin{example}\label{ex:quadrangle-ambiguity}\mbox{}
  \begin{enumerate}[(a)]
  \item\label{item:23} Fix a trapezoidal face $\bigtriangleup$ of $\Gamma$ with vertices as in
    \ref{lem:1}. One can remove the vertex $D$ if and only if either $n=1$ or
    $r_2+tq=1$.  The vertex \(A\) can be removed if and only if $r_2+q=1$.
    (There are analogous characterizations for $B$ and $C$, too.) The case
    $n=1$ is the `moving triangle' situation \ref{ex:moving-triangle}. If
    $r_2+q=1$, then there are (at least) two possibilities for the choice of
    the axis of $\kn$, namely the segments $[(0,1,n),(r_1+tp,0,0)]$ and
    $[(p,0,n),(r_1+tp-p,1,0)]$. One of them replaces the trapezoid by a
    triangle, while the other replaces it by a smaller trapezoid. Hence, in
    any situation, if a trapezoid can be decreased in some way, then it can be
    replaced by a triangle in the equivalence class of the diagram. Otherwise,
    it is called \emph{non-removable} (this happens if $n>1$, $r_1+p>1$,
    $r_2+q>1$).
  \item\label{item:24} If above $q=p=1$ and $r_1=r_2=0$, then $\bigtriangleup$ is the
    canonical representative of its class. One has four possible axes, and $\bigtriangleup$
    can be reduced to the trapezoid $(0,1,n),(1,0,n),(t-1,1,0),(1,t-1,0)$ or
    to the triangles $(0,1,n),(t,0,0),(1,t-1,0)$ or
    $(1,0,n),(t-1,1,0),(0,t,0)$. These are the minimal representatives.
  \end{enumerate}
\end{example}

\begin{remark}\label{re:99}\mbox{}
  \begin{enumerate}[(a)]
  \item\label{item:25} Let $\bigtriangleup$ be a triangular face of an M1-minimal
    representative $\Gamma$. Then in any minimal representative of $\Gamma$, which is
    obtained from $\Gamma$ via moves $\kn$, $\bigtriangleup$ survives as a triangular face which
    is independent of the choice of the minimal representative.  This happens,
    because the axes of all the moves $\kn$, which can be applied to $\bigtriangleup$,
    cannot intersect each other, hence all of them can be applied
    `simultaneously' (a fact, which is not true in the case of removable
    trapezoids, see \ref{ex:quadrangle-ambiguity} above).

  \item\label{item:26} Therefore, any class whose \emph{canonical}
    representative has a non-removable trapezoid, or a central triangle or a
    central edge, admits a unique minimal representative
  \end{enumerate}
\end{remark}

\begin{disdef}[d-minimal representatives]\label{cor:str2} Fix a class.
  It may contain many minimal representatives; we will distinguish one of
  them, and we call it \emph{d-minimal} (distinguished-minimal). If the class
  admits a unique minimal representative, then there is no ambiguity for the
  choice. This happens e.g.\ in all the situations \ref{re:99}(\ref{item:26}).

  For the sake of completeness, we allow diagrams satisfying \ref{cor:3}.  For
  such a class, the (d-)minimal representative is the segment
  $[(0,1,1),(n,0,0)]$, for some $n\geq 2$, as given in
  \ref{cor:3}(\ref{item:15}).

  Next, assume that a canonical representative contains a \emph{removable
    trapezoid} (i.e.\ one replaceable by a triangle).  Using the notations
  of~\ref{lem:1}, if $n>1$, then again there is a unique minimal
  representative, unless we are in the situation of
  \ref{ex:quadrangle-ambiguity}(\ref{item:24}) (when there are two, but they
  correspond to each other by a permutation of coordinates). By definition,
  this is the d-minimal representative (in the last case it is well-defined up
  to the permutation of coordinates).

  If $n=1$, then we are in the situation of a moving triangle
  \ref{ex:moving-triangle}, and the class may contain many minimal
  representatives. (An even more annoying fact is that such a class may
  contain two equivalent diagrams such that one of them has a central triangle
  while the other has a central edge.) We will declare the position of the
  moving point $R$ for the d-minimal representative as follows. Assume that
  $p<q$ (for $q<p$ interchange $z_1$ and $z_2$). If $R$ cannot be moved to any
  of the coordinate axis (cf.~\eqref{eq:movcor}), then take for $R$ that
  possible lattice point which is closest to the \(z\sb{1}\) axis. If $R$ can
  be moved to exactly one coordinate axis, then move it there.  If $R$ can be
  moved to both axes, then move to the \(z\sb{1}\) axis.  (Since the
  determinants of $QR$ and $PR$ are $p$ and $q$, respectively, by this choice
  of $R$, the determinant of the edge lying on the coordinate plane is larger.
  There is no deep motivation for this choice, except that we need one.  In
  the `inverse' algorithm the very same choice is built in.)
\end{disdef}

\begin{corollary}[Structure of d-minimal representatives]\label{cor:str}
  (Cf.\ also with \ref{prop:1}.) The d-minimal representative of an
  equivalence class (which does not satisfy \ref{cor:3}) sits in exactly one
  of the following three disjoint families of diagrams characterized by the
  existence of
\begin{enumerate}[(1)]
\item\label{item:27} a non-removable trapezoid,
\item\label{item:28} a central triangle, or
\item\label{item:29} a central edge.
\end{enumerate}
\end{corollary}

\begin{notation}\label{symbols}
  The three disjoint families listed in \ref{cor:str} will be denoted by
  $\nrt$, $\ctr$, $\ced$. They can be divided further according to the number
  of hands. This number will appear as a subscript. E.g., $\nrt_3$ denotes
  that family of classes of Newton boundaries whose d-minimal representative
  has a non-removable trapezoid and \(3\) hands.
\end{notation}

The first fruit of the minimality of a graph is the following arithmetical
criterion:
\begin{proposition}\label{negy}
  Fix a minimal representative $\Gamma_{\min}$ of a class which does not satisfy
  \ref{cor:3}. Consider an edge of $\partial\Gamma_{\min}$ which is the intersection of
  the faces $\bigtriangleup$ and $\bigtriangledown$ of $\Gamma_{\min,+}$, where the second one is non-compact.
  Then $\innerdeterminant{\bigtriangleup}{\bigtriangledown}>1$.
\end{proposition}
\begin{proof}
We have to analyze two types of edges, cf.~\ref{rem:uj}. First we discuss
edges on a coordinate plane, say $[(q_1,0,q_3),(q_1',0,q_3')]$.  Take a
triangle in $\bigtriangleup$ which satisfies the criterions of \ref{lem:9}(\ref{item:79}).
Then $\innerdeterminant{\bigtriangleup}{\bigtriangledown}=
\innerdeterminant{\bigtriangleup}{\coordinatevector{2}}=p_2$. But, if $p_2=1$, then this
triangle can be eliminated by $\jn$.

Now we turn to the other type of edges, which have the form
\(AB=[(a,0,c),(0,1,b)]\) with $a>0$. Take a third vertex $C=(r,s,u)$ on $\bigtriangleup$
such that \(\trianglebyvertices{A}{B}{C}\) is empty. Then the identity
\eqref{eq:24} of \ref{lem:10} can be applied:
$\innerdeterminant{\bigtriangleup}{\bigtriangledown}=r+(s-1)a$.  Assume that $r+(s-1)a=1$.  If $s=0$ then
$[(r,0,u),(0,1,b)]$ is an axis of a move $\jn$, hence
\(\trianglebyvertices{A}{B}{C}\) can be eliminated.  If $s=1$, then $r=1$,
hence by \eqref{eq:rhs} $u=0$, which contradicts our assumption about
\ref{cor:3}.  The remaining case $s\geq 2$ imposes $r=0$, $s=2$, $a=1$, with
$[(1,0,c),(0,2,u)]$ an axis of $\jn$ which eliminates
\(\trianglebyvertices{A}{B}{C}\).
\end{proof}

\begin{remark}\label{rem:200}\mbox{}
  \begin{enumerate}[(a)]
  \item\label{item:30} By the above proof, when we eliminate triangles from a
    diagram by moves $\jn$, then, in fact, we eliminate those `mixed
    determinants' (i.e.\ when a face is non-compact) with
    $\innerdeterminant{\bigtriangleup}{\bigtriangledown}=1$. By \ref{negy}, by repeated application of
    $\jn$, we can eliminate all such mixed determinants, provided that the
    class does not satisfies \ref{cor:3}.  (Otherwise this is not true:
    \ref{ex:105} shows a minimal triangle with a `mixed determinant' \(1\).)
  \item\label{item:31} The statement of \ref{negy} is also true for a class
    which satisfies \ref{cor:3} (where $\bigtriangleup$ and $\bigtriangledown$ are non-compact and contain
    $[(0,1,1),(n,0,0,)]$): $\innerdeterminant{\bigtriangleup}{\bigtriangledown}=n>1$.
  \end{enumerate}
\end{remark}

\section{The dual resolution graph}
\label{sec:dual-resol-graph}

\subsection{Graph terminology}
\label{sec:graph-terminology}

\begin{para}\label{gr10}
  Recall that any resolution graph $G(f)$ of $\{f=0\}$ is also a possible
  plumbing graph of the link $K(f)$ of $f$. The link $K(f)$ is a rational
  homology sphere if and only if $G(f)$ is a tree, and the genera of all the
  vertices are \(0\). In such a case, $G(f)$ has only one set of decorations:
  each vertex carries the self-intersection number of the corresponding
  irreducible exceptional divisor. In this subsection we recall the
  terminology of resolution graphs, and we present a construction which
  `simplifies' a given graph. Its output will be called the \emph{orbifold
    diagram}.
\end{para}

\begin{para}\label{gr11}
  Let $G$ be a decorated tree with vertices $\calv$ and decorations
  ${\{b_v\}}_{v\in\calv}$. The entries of \emph{intersection matrix}
  $(I_{vw})_{v,w\in\calv}$ of $G$ are $I_{vv}=b_v$, and for $v\neq w$ one sets
  $I_{vw}=1$ if $[vw]$ is an edge, and $I_{vw}=0$ otherwise. We assume that
  $I$ is negative definite (since the matrix of a dual resolution graph is so
  \cite{Mumford}). By definition, $\det(G)\coloneqq \det(-I)$ is the
  \emph{determinant of the graph} $G$.

  A \emph{node} of $G$ is a vertex whose degree is at least \(3\).  Let
  $\calr$ be their collection. A \emph{chain} is the path between two nodes
  excluding the endpoints, which does not contain any nodes.  (We say that the
  chain \emph{connects} the two nodes.) Similarly, a \emph{leg} of $G$ is a
  path between a degree \(1\) vertex and a node containing the degree \(1\)
  vertex but not the node, and containing no other nodes, either.

  If $r,s\in\calr$ are connected by a chain in $G$, then the determinant of this
  chain (i.e., the determinant of the corresponding subgraph) will be denoted
  by $n_{rs}$.

  A \emph{star-shaped} graph is a graph with a unique node. For any $r\in\calr$
  there is a unique maximal star-shaped subgraph $G_r$ of $G$ which contains
  $r$.

  In general, a star-shaped graph is a plumbing graph of a Seifert
  \(3\)-manifold. This has a natural $S^1$-action and orbifold structure. If
  the star-shaped graph $G_r$ has normalized Seifert invariants, say,
  ${(\alpha_i,\omega_i)}_i$ (here, each pair is associated with one of the legs of the
  subgraph, $\alpha_i$ is the leg-determinant, $0\leq \omega_i<\alpha_i$, and we put the pair
  $(1,0)$ for legs with determinant one), and central vertex with decoration
  $b_r$, then the \emph{orbifold Euler number} of $G_r$ is $e_r\coloneqq
  b_r+\sum_i\omega_i/\alpha_i$, see e.g.~\cite{Wa} for details. \end{para}

\begin{para}[The orbifold diagram]\label{gr12}
  Sometimes we do not need all the data of $G$, but only its shape and the
  determinants of some of its subgraphs. This information will be codified in
  a simpler graph-like diagram, the \emph{orbifold diagram} associated with
  $G$, denoted by $G^o$.

  $G^o$ is constructed from $G$ as follows.  $G^o$ has vertices, edges
  connecting two vertices, and half-free edges.  A half-free edge is attached
  with one of its ends to a vertex, while its other end is free. The vertices
  of $G^o$ are the nodes of \(G\).  The (ordinary) edges of \(G^o\) are the
  chains of \(G\).  The endpoints of an edge are the two nodes it connects as
  a chain.  The half-free edges are the legs of \(G\).  The endpoint of a
  half-free edge is the node to which it is adjacent in \(G\) as a leg.  Then
  we decorate $G^o$: we put on each edge the determinant of the corresponding
  chain or leg, and we label each node $r$ with the orbifold Euler number
  $e_r$ of the star-shaped subgraph $G_r$. (In the special case when $G$ has
  no nodes then $G^o$ is a `free' edge decorated by $\det(G)$.)

  The half-free edges of the orbifold diagram $G^o$ will still be called legs.

  The entries of the \emph{orbifold intersection matrix}
  $(I^o_{rs})_{r,s\in\calr}$ of $G^o$, by definition, are $I^o_{rr}=e_r$, and
  for $r\neq s$ one sets $I^o_{rs}=1/n_{rs}$ if $[rs]$ is an edge of $G^o$, and
  $I^o_{rs}=0$ otherwise. (Here we will not explain the `orbifold geometry'
  behind this definition.  Nevertheless, for a possible motivation, see
  §\ref{sec:splice-diagram}.) Similarly as above, we set $\det(G^o)\coloneqq
  \det(-I^o)$.
\end{para}
\begin{lemma}\label{gr13}
  Fix a graph $G$ as above with $\calr\neq\emptyset$. Let $\Pi$ be the product of the
  determinants of all the chains and legs of $G$. Then $I^o$ is negative
  definite, and
  \begin{equation}\label{eq:gr}
    \det(G)=\det(G^o)\cdot\Pi.
  \end{equation}
\end{lemma}
\begin{proof}
The negative definiteness of $I^o$ follows from \eqref{eq:gr} applied to some
subgraphs.  The equality~\eqref{eq:gr} is elementary linear algebra, it
follows (e.g.) by induction on $\#\calr$. If $\#\calr=1$, then \eqref{eq:gr} is
well-known, see e.g.\ \cite{Neumann}. The induction runs as follows. Fix
$r,s\in\calr$ which are connected by a chain $G_{rs}$. The connected components
of $G\setminus (\{r,s\}\cup G_{rs})$ are ${\{G_i\}}_i$, the connected component of $G\setminus \{r\}$
which contains $s$ is $G_{(s)}$, and similarly one defines $G_{(r)}$. Then
$\det(G)\cdot \det(G_{rs})=\det(G_{(s)})\cdot\det(G_{(r)})-\prod_i\det(G_i)$.
\end{proof}
\begin{remark}\label{NW}
  The orbifold diagram has exactly the same shape as the \emph{splice diagram}
  considered in \cite{NW}, but it has different decorations. Nevertheless, by
  similar identities what we used in the proof of \ref{gr13}, one can show
  that the orbifold diagram contains the same amount of information as the
  splice diagram and $\det(G)$ altogether.
\end{remark}

\subsection{Oka's algorithm for $G(f)$. The case of minimal
  representatives.}\label{OKAsect}

\begin{para}
  Let $f\colon (\setC^3,0) \to (\setC,0)$ be a germ with \emph{isolated singularity} and
  \emph{non-degenerate Newton principal part} whose link $K(f)$ is a rational
  homology sphere. In particular, its Newton boundary $\Gamma(f)$ satisfies
  \eqref{eq:2} and \eqref{eq:rhs}. In the first part of this subsection we
  recall the combinatorial algorithm of M. Oka \cite[Theorem~6.1]{MR894303},
  which provides a (possible, in general non-minimal) dual resolution graph
  $G(f)$ of the surface singularity $(\{f=0\},0)$ from $\Gamma(f)$.

  In order to emphasize the dependence of the output upon $\Gamma(f)$, we write
  $G(\Gamma(f))$.
\end{para}
\begin{para}[Notations]
  Recall that $\Gamma(f) $ is the union of compact faces of $\Gamma_+\coloneqq
  \Gamma_+(\supp(f))$, which can be recovered from $\Gamma(f)$ as $\Gamma_+(\{\text{vertices
    of $\Gamma(f)$}\})$.  Hence, they contain the same amount of information.
  Similarly as above, $\calf$ denotes the collection of all faces of $\Gamma_+$,
  and $\calf_c$ denotes the set of all compact faces of $\Gamma_+$.  For any
  $\bigtriangleup\in\calf_c$, we write $\calf_\bigtriangleup$ for the collection of all faces of $\Gamma_+$
  adjacent to $\bigtriangleup$.  Other notations are from~\ref{not}.
\end{para}
\begin{para}[The algorithm]\label{alg}
  The graph $G(\Gamma(f))$ is a subgraph of a larger graph $\widetilde{G}(\Gamma(f))$,
  whose construction is the following. To start with, we consider $\calf$ as a
  set of vertices (we will call them \emph{face vertices}).  Then, if
  $\bigtriangleup,\bigtriangledown\in\calf$ are two adjacent faces, then we connect them by
  $\multiplicity{\bigtriangleup}{\bigtriangledown}$ copies of the following chain.

  If \(\innerdeterminant{\bigtriangleup}{\bigtriangledown}>1\) then let \(0 < \subchaindeterminant{\bigtriangleup}{\bigtriangledown} <
  \innerdeterminant{\bigtriangleup}{\bigtriangledown}\) be the unique integer for which
  \begin{equation}\label{eq:5} {\namedvector{c}}_{\bigtriangleup,\bigtriangledown} \coloneqq
    (\normalvector{\bigtriangledown} +
    \subchaindeterminant{\bigtriangleup}{\bigtriangledown}
    \normalvector{\bigtriangleup})/\innerdeterminant{\bigtriangleup}{\bigtriangledown}
  \end{equation}
  is an integral vector.  Let us write \(\innerdeterminant{\bigtriangleup}{\bigtriangledown}/
  \subchaindeterminant{\bigtriangleup}{\bigtriangledown}\) as a continued fraction:
  \begin{equation}
    \label{eq:6}
    \frac{\innerdeterminant{\bigtriangleup}{\bigtriangledown}}
    {\subchaindeterminant{\bigtriangleup}{\bigtriangledown}} =
    b_1 - \cfrac{1}{b_2 -\cfrac{1}{ \dotsb - \cfrac{1}{b_k}}},
  \end{equation}
  where each \(b_i \geq 2\).  Then the chain with the corresponding
  self-intersection numbers is

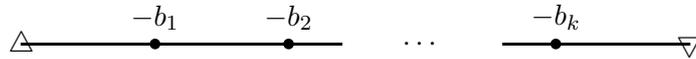
\begin{figure}[h]
  \centering
  \begin{picture}(270,30)(90,5) \put(150,20){\circle*{4}}
    \put(200,20){\circle*{4}} \put(300,20){\circle*{4}}
    \put(100,20){\line(1,0){120}} \put(280,20){\line(1,0){70}}
    \put(100,20){\makebox(0,0){$\bigtriangleup$}} \put(350,20){\makebox(0,0){$\bigtriangledown$}}
    \put(150,30){\makebox(0,0){$-b_1$}} \put(200,30){\makebox(0,0){$-b_2$}}
    \put(300,30){\makebox(0,0){$-b_k$}} \put(250,20){\makebox(0,0){$\cdots$}}
  \end{picture}
  \caption{Chain between two face vertices}
  \label{eq:pict}
\end{figure}
The left ends of all the $\multiplicity{\bigtriangleup}{\bigtriangledown}$ copies of the chain (marked by
$\bigtriangleup$) are identified with the face vertex corresponding to $\bigtriangleup$, and similarly
for the right ends marked by $\bigtriangledown$.

If \(\innerdeterminant{\bigtriangleup}{\bigtriangledown}=1\) then the chain consists of an edge connecting
the vertices \(\bigtriangleup\) and \(\bigtriangledown\) (we put $\multiplicity{\bigtriangleup}{\bigtriangledown}$ of them). Also, in
this case we set $\subchaindeterminant{\bigtriangleup}{\bigtriangledown} \coloneqq 0$ and
${\namedvector{c}}_{\bigtriangleup,\bigtriangledown} \coloneqq \normalvector{\bigtriangledown}$.

Next, we compute the decoration \(\selfintersectionnumber{\bigtriangleup}\) of any face
vertex \(\bigtriangleup\in\calf_c\) by the equation:
\begin{equation}
  \label{eq:7}
  \selfintersectionnumber{\bigtriangleup} \normalvector{\bigtriangleup} +
  \sum_{\bigtriangledown\in\calf_\bigtriangleup}
  \multiplicity{\bigtriangleup}{\bigtriangledown}
  {\namedvector{c}}_{\bigtriangleup,\bigtriangledown} = \namedvector{0}.
\end{equation}
What we get in this way is the graph $\widetilde{G}(\Gamma(f))$.  Notice that the
face vertices corresponding to \emph{non-compact faces} are not decorated. If
we delete all these vertices (and all the edges adjacent to them) we get the
\emph{dual resolution graph $G(\Gamma(f))$}.

Notice that $\widetilde{G}(\Gamma(f))$ has no legs, but some of the chains of
$\widetilde{G}(\Gamma(f))$ become legs of $G(\Gamma(f))$ when we delete the vertices
corresponding to $\calf\setminus \calf_c$.  Regardless whether a chain in
Figure~\ref{eq:pict} transforms into a chain or a leg of $G(\Gamma(f))$, it keeps
its determinant $\innerdeterminant{\bigtriangleup}{\bigtriangledown}$.
\end{para}
\begin{remark}
  If one starts with another Newton diagram, say $\Gamma'(f)$, obtained from $\Gamma(f)$
  via Moves~\ref{item:6} or~\ref{item:7}, then the graph $G(\Gamma'(f))$ can be
  obtained from $G(\Gamma(f))$ by some blow-ups, in accordance with
  \ref{2.8}(\ref{item:9}). Hence, in general, $G(\Gamma(f))$ is not a good minimal
  resolution graph.  Recall that a dual resolution graph $G(f)$ with all
  genera vanishing is good minimal if all its $(-1)$-vertices are nodes.  Each
  normal surface singularity admits a \emph{unique good minimal
    resolution}.\end{remark}
\begin{proposition}\label{mingood}
  If the Newton diagram $\Gamma_{\min}(f)$ is a minimal representative of its class
  then the output $G(\Gamma_{\min}(f))$ of Oka's algorithm is the good minimal
  resolution graph. In fact, $G(\Gamma_{\min}(f))$ reflects the shape of the
  diagram $\Gamma_{\min}(f)$ (preserving the corresponding adjacency relations):
  \begin{enumerate}[(a)]
  \item\label{item:32} the nodes of $G(\Gamma_{\min}(f))$ correspond bijectively to
    the faces of $\Gamma_{\min}(f)$;
  \item\label{item:33} the chains and legs of $G(\Gamma_{\min}(f))$ correspond
    bijectively to the edges of $\Gamma_{\min}(f)$ \emph{not lying} in $\partial
    \Gamma_{\min}(f)$ and the edges \emph{lying} on $\partial \Gamma_{\min}(f)$, respectively.
    (In the case of~\ref{cor:3} we understand by this that $\Gamma=\partial\Gamma$ is a
    segment, and $G(\Gamma_{\min}(f))$ is a string.)
  \end{enumerate}
\end{proposition}
\begin{proof}
The chains in Figure~\ref{eq:pict} contain no $(-1)$-vertex, any face has at
least three edges, and all the leg-determinants are greater than \(1\) by
\ref{negy} and \ref{rem:200}, hence the statement follows.
\end{proof}
The legs corresponding to different primitive segments of the \emph{same} edge
form a \emph{leg group}.

In the next subsection we make a more direct connection between the normal
vectors of faces, the coordinates of vertices of \(\Gamma_{\min}(f)\), and the
determinants of legs in $G(\Gamma_{\min}(f))$.

\subsection{Leg-determinants in $G(\Gamma_{\min}(f))$} \label{sec:determinant-leg}

We fix a minimal Newton diagram $\Gamma_{\min}=\Gamma_{\min}(f)$ which does not satisfy
\ref{cor:3}, and let $\bigtriangleup\in\calf_c$ be one of its faces.  Let us consider the
legs in $G(\Gamma_{\min})$ adjacent to $\bigtriangleup$.  By~\ref{mingood} they correspond to
the primitive segments lying on the edges of $\partial \Gamma_{\min}\cap\bigtriangleup$. The next
proposition summarizes the divisibility properties of the determinants of
these legs. We will refer to such a leg-determinant as the \emph{determinant}
$D(\alpha)$ of the corresponding edge $\alpha$. (The coordinate choices are accidental,
they can be permuted arbitrarily.)  As usual, we write
\(\normalvector{\bigtriangleup}=(a_1,a_2,a_3)\).

\begin{proposition} \label{prop:2} \mbox{}
  \begin{enumerate}[I.]
  \item\label{item:34} Any edge $\alpha$ of $\partial \Gamma_{\min}\cap\bigtriangleup$ satisfies the next
    divisibility properties:
    \begin{enumerate}[1.]
    \item\label{item:35} If $\alpha\subset \{z_3=0\}$, then $D(\alpha)=\mygcd{a_1,a_2}$ and
      $\mygcd{D(\alpha),a_3}=1$.
    \item\label{item:36} If \(\alpha=[(p,0,a),(0,1,b)]\) crosses the \(z\sb{3}\)
      axis, then $D(\alpha)=a_3$, and $D(\alpha)$ does not divide any of $\{a_1,a_2\}$,
      unless the edge can be `moved' to a coordinate plane,
      see~\ref{ex:moving-triangle}.  If this happens, then $D(\alpha)$ divides the
      corresponding two coordinates as in part~(\ref{item:35}), and either of
      the following cases holds.
      \begin{enumerate}[(a)]
      \item\label{item:37} The edge belongs to a `moving triangle' so that it
        can be moved to a coordinate plane by moving the moving vertex to a
        coordinate axis, see~\eqref{eq:movcor}.
      \item\label{item:38} \((0,2,c)\) is on \(\bigtriangleup\), hence \(\bigtriangleup\) can be
        extended (in the class of $\Gamma_{\min}$) with a new vertex
        \((0,0,2b-c)\), which lengthens the edge $\beta= [(0,1,b),(0,2,c)]$ by an
        extra primitive segment.  We interpret this as moving $\alpha$ to this
        extra segment whose determinant is $D(\beta)$.
      \end{enumerate}
    \end{enumerate}
  \item\label{item:39} The determinants belonging to different edges of \(\bigtriangleup\)
    \begin{enumerate}[1.]
    \item\label{item:40} differ, except in the case of~(\ref{item:38}), where
      $D(\alpha)=D(\beta)$;
    \item\label{item:41} are pairwise relative prime except the two cases
      below:
      \begin{enumerate}[(a)]
      \item\label{item:42} an edge lying on a coordinate plane is adjacent to
        a crossing edge: the determinant of the former one divides the
        determinant of the latter one;
      \item\label{item:43} an edge $\alpha$ lying on a coordinate plane is adjacent
        to two crossing edges: then $D(\alpha)$ is the greatest common divisor of
        the determinants of the crossing edges.
      \end{enumerate}
    \end{enumerate}
  \end{enumerate}
\end{proposition}

\begin{proof} (\ref{item:35}) is clear by \eqref{eq:ijk} since
\(\normalvector{\bigtriangleup}\) is primitive.  \eqref{eq:23} implies the first part of
(\ref{item:36}). Since \(\normalvector{\bigtriangleup}\) is orthogonal to $\alpha$, i.e.\
$\scalarproduct{\normalvector{\bigtriangleup}}{(p,0,a)}
=\scalarproduct{\normalvector{\bigtriangleup}}{(0,1,b)}$, one has
\begin{equation}
  \label{eq:8}
  a_{2} = p a_{1} + (a-b) a_{3}.
\end{equation}
Since \(\normalvector{\bigtriangleup}\) is primitive $a_3\nmid a_1$, and if $a_3\mid a_2$ then
$D(\alpha)=a_3\mid p$, too. This will be used later.

Let \(P\) be a lattice point of \(\bigtriangleup\), such that the triangle \(\bigtriangledown\) (a part of
the face \(\bigtriangleup\)) formed by $\alpha$ and \(P\) is empty.  First notice that \(P\)
cannot be on the \(z\sb{1} z\sb{3}\) plane, since then \(\bigtriangledown\) would be
removable. If \(P\) lies on the \(z\sb{2} z\sb{3}\) plane and has the form
\((0,q,c)\), then by \eqref{eq:24}, $D(\alpha)=(q-1)p$. If $D(\alpha)\mid p$ then $q=2$.
This is the case~(\ref{item:38}).

Assume now that \(P\) lies on the \(z\sb{1} z\sb{2}\) plane, \(P=(r_1,r_2,0)\)
with \(r_1\) and \(r_2\) positive.  Then, again by \eqref{eq:24}, $D(\alpha)=r_1+ p
(r_2 - 1)$. Thus, if $D(\alpha)\mid p$ then \(r_2=1\) and \(r_1\mid p\).  But then the
triangle \(\bigtriangledown\) is movable (as in~\ref{ex:moving-triangle}), and the vertex
\((p,0,a)\) can be moved to the point \((0, 0, a + b p/r_1)\) lying on the
third coordinate axis.  This is the case~(\ref{item:37}).

For part~(\ref{item:39}), assume that the leg-determinants belonging to two
different edges are not relative prime.  If one of the edges lie on a
coordinate plane, then \eqref{eq:ijk} and \eqref{eq:23} show that we are in
the situation~(\ref{item:42}) or~(\ref{item:43}).  Otherwise, if one of the
edges crosses, say, the \(z\sb{1}\) axis, and the other edge crosses the
\(z\sb{2}\) axis, then their endpoints sitting on the \(z\sb{1} z\sb{2}\)
plane do not coincide, and hence case~(\ref{item:43}) holds. Indeed, assume
that the two endpoints do coincide. This common point cannot be \((1,1,0)\) by
our assumption, cf.~\ref{cor:3}. Otherwise, by a relation similar to
\eqref{eq:8}, the third coordinate of \(\normalvector{\bigtriangleup}\) is an integral
linear combination of the first two ones (which are the determinants),
contradicting the fact that \(\normalvector{\bigtriangleup}\) is primitive.

Finally, if two determinants are equal, then by part~(\ref{item:41}) the
corresponding edges must be as in (\ref{item:42}). But then, by \eqref{eq:p2}
and \eqref{eq:24}, we have $q=2$ which leads to~(\ref{item:38}).
\end{proof}

\begin{corollary}[Non-removable trapezoids]
  \label{lem:4} The leg groups of a non-removable trapezoid have different
  determinants. Hence, the collection of chains and legs adjacent to the
  vertex corresponding to the trapezoid can be separated in \(4\)
  distinguishable groups.
\end{corollary}
Since a vertex corresponding to a triangular face has at most \(3\) such
groups of distinguishable legs and chains, the vertex of a non-removable
trapezoid can be recognized in the resolution graph.
\begin{proof}[Proof of \ref{lem:4}] Assume the contrary. Then by~\ref{prop:2},
we are in the situation of~(\ref{item:38}) with the points $(p,0,a)$,
$(0,1,b)$ and $(0,2,c)$ on the trapezoid (modulo a permutation of the
coordinates). Then this face can be extended by the vertex $(0,0,2b-c)$. This
extended face is a trapezoid, too. But
by~\ref{ex:quadrangle-ambiguity}(\ref{item:23}), this is a removable
trapezoid.
\end{proof}

\subsection{The orbifold diagram} \label{sec:splice-diagram}

\begin{para}\label{or10}
  \emph{We fix a minimal representative} $\Gamma_{\min}(f)$ as in \ref{mingood}.
  §\ref{OKAsect} provides a good minimal resolution graph $G(\Gamma_{\min}(f))$
  from the Newton diagram $\Gamma_{\min}(f)$. On the other hand, to any graph $G$,
  the general procedure \ref{gr12} associates a diagram $G^o$. In the present
  situation this will be denoted by $G^o(\Gamma_{\min}(f))$. Although, by the very
  construction of $G^o$, we (apparently) throw away some information, we
  prefer to use $G^o(\Gamma_{\min}(f))$ since it reflects more faithfully the
  Newton diagram. For the convenience of the reader, in short, we sketch how
  one can draw $G^o(\Gamma_{\min}(f))$ directly from $\Gamma_{\min}(f)$.

  Similarly as in §\ref{OKAsect}, first we construct a decorated graph
  $\widetilde{G}^o$.  Its vertices are the elements of $\calf$, i.e.\ all the
  faces of $\Gamma_{\min,+}$.  If \(\bigtriangleup, \bigtriangledown\in\calf\) are adjacent in $\Gamma_{\min,+}$, then
  we connect them by \(\multiplicity{\bigtriangleup}{\bigtriangledown}\) edges in $\widetilde{G}^o$, and
  we label each of these edges with the number \(\innerdeterminant{\bigtriangleup}{\bigtriangledown}\).
  Finally, we label each \(\bigtriangleup\in\calf_c\) with the orbifold Euler number of the
  maximal star-shaped subgraph containing \(\bigtriangleup\), which is
  \begin{equation}
    \label{eq:9}
    \Eulernumber{\bigtriangleup} \coloneqq \selfintersectionnumber{\bigtriangleup}
    +\sum_{\bigtriangledown\in\calf_\bigtriangleup}
    \multiplicity{\bigtriangleup}{\bigtriangledown}
    \frac{\subchaindeterminant{\bigtriangleup}{\bigtriangledown}}{\innerdeterminant{\bigtriangleup}{\bigtriangledown}}.
  \end{equation}
  In this way we get the labelled graph $\widetilde{G}^o$. If we remove the
  vertices \(\{v_\bigtriangleup : \bigtriangleup\in\calf\setminus \calf_c\}\) (but we keep the edges---i.e.\ the new
  legs---adjacent to them), we get the diagram $G^o(\Gamma_{\min}(f))$.  For any
  $\bigtriangleup\in\calf_c$, we call $\Eulernumber{\bigtriangleup}$ the \emph{orbifold Euler number of}
  $\bigtriangleup$.
\end{para}

\begin{para}
  The point is that \eqref{eq:7} can be transformed via the orbifold Euler
  numbers into some (more natural) identities which only involve the normal
  vectors of the faces.
\end{para}

\begin{proposition}\label{lem:5} Fix the representative $\Gamma_{\min}(f)$. Then
  for any $\bigtriangleup\in\calf_c$ one has:
  \begin{equation}
    \label{eq:10}
    \Eulernumber{\bigtriangleup} \normalvector{\bigtriangleup}
    +\sum_{\bigtriangledown\in\calf_\bigtriangleup}
    \frac{\multiplicity{\bigtriangleup}{\bigtriangledown}}{\innerdeterminant{\bigtriangleup}{\bigtriangledown}}
    \normalvector{\bigtriangledown} = \namedvector{0}.
  \end{equation}
\end{proposition}
\begin{proof} Use \eqref{eq:5}, \eqref{eq:7} and \eqref{eq:9}.
\end{proof}

Obviously, if one wishes to recover the equation of a face of $\Gamma_{\min}(f)$,
one needs its normal vector \(\normalvector{\bigtriangleup}\), and its \emph{face value},
i.e.\ the value of \(\normalvector{\bigtriangleup}\) on any of the face's point:
\begin{equation}
  \label{eq:11}
  \facevalue{\bigtriangleup} \coloneqq
  \scalarproduct{\normalvector{\bigtriangleup}}{P}, \quad P \in \bigtriangleup.
\end{equation}
It turns out that these numbers ${\{\facevalue{\bigtriangleup}\}}_{\bigtriangleup\in \calf}$ also satisfy a
similar equation:
\begin{proposition}\label{lem:6}
  Fix $\Gamma_{\min}(f)$ as above. Then for any $\bigtriangleup\in\calf_c$ one has:
  \begin{equation}
    \label{eq:12}
    \Eulernumber{\bigtriangleup} \facevalue{\bigtriangleup} +
    \sum_{\bigtriangledown\in\calf_\bigtriangleup}
    \frac{\multiplicity{\bigtriangleup}{\bigtriangledown}}{\innerdeterminant{\bigtriangleup}{\bigtriangledown}} \facevalue{\bigtriangledown} =
    - \combinatorialarea{\bigtriangleup}.
  \end{equation}
\end{proposition}

\begin{proof} Denote the vertices of \(\bigtriangleup\) by \(P_0, \dotsc, P_k\) in this
order, and set \(P_{k+1} \coloneqq P_0\). Assume that $P_iP_{i+1}$ is the
common edge of \(\bigtriangleup\) and \(\bigtriangledown\).  Then, by \eqref{eq:21}, one has
\begin{equation}\label{eq:13}
  \combinatorialarea{\trianglebyvertices{P_0}{P_{i}}{P_{i+1}}} =
  \frac{\multiplicity{\bigtriangleup}{\bigtriangledown}}{\innerdeterminant{\bigtriangleup}{\bigtriangledown}}
  \scalarproduct{\normalvector{\bigtriangledown}}{P_0-P_i}=
  \frac{\multiplicity{\bigtriangleup}{\bigtriangledown}}{\innerdeterminant{\bigtriangleup}{\bigtriangledown}}
  \scalarproduct{\normalvector{\bigtriangledown}}{P_0}-
  \frac{\multiplicity{\bigtriangleup}{\bigtriangledown}}{\innerdeterminant{\bigtriangleup}{\bigtriangledown}}
  \facevalue{\bigtriangledown},
\end{equation}
with \(\combinatorialarea{\trianglebyvertices{P_0}{P_{0}}{P_{1}}}=
\combinatorialarea{\trianglebyvertices{P_0}{P_{k}}{P_{k+1}}}=0\). Then, by
\eqref{eq:10} and \eqref{eq:13}
\begin{multline*}
  -\Eulernumber{\bigtriangleup} \facevalue{\bigtriangleup} = -\Eulernumber{\bigtriangleup}
  \scalarproduct{\normalvector{\bigtriangleup}}{P_0}= \\
  \sum \frac{\multiplicity{\bigtriangleup}{\bigtriangledown}}{\innerdeterminant{\bigtriangleup}{\bigtriangledown}}
  \scalarproduct{\normalvector{\bigtriangledown}}{P_0}= \sum \frac{\multiplicity{\bigtriangleup}{\bigtriangledown}}
  {\innerdeterminant{\bigtriangleup}{\bigtriangledown}} \facevalue{\bigtriangledown}+\sum
  \combinatorialarea{\trianglebyvertices{P_0}{P_{i}}{P_{i+1}}}.
\end{multline*}
Then use the additivity of the combinatorial area.
\end{proof}
\begin{corollary}\label{cor:negdef}
  For any subset ${\overline{\calf}}_c\subset \calf_c$ consider the system of
  equations \eqref{eq:10} for all $\bigtriangleup\in{\overline{\calf}}_c$ in unknowns
  $\{\normalvector{\bigtriangleup} : \bigtriangleup\in{\overline{\calf}}_c\}$ (or in one of their fixed
  coordinates). Then this system is non-degenerate. The same is true for
  equations \eqref{eq:12} instead of \eqref{eq:10}.
\end{corollary}
\begin{proof}
The matrix of the system (for ${\overline{\calf}}_c=\calf_c$) coincides with
the matrix \(I^o\) of the orbifold diagram \(G^o(\Gamma_{\min}(f))\)
(cf.~\ref{gr12}), which is negative definite by~\ref{gr13}. \end{proof}

\begin{remark}\label{lem:11b} If one wishes to solve the above equations, one needs
  the values for non-compact faces. If such a face is supported by a
  coordinate plane, then its normal vector is a coordinate vector, and its
  face value is \(0\). Otherwise, if it has an edge of type
  \([(a,0,c),(0,1,b)]\) ($a>0$), then its normal vector is $(1,a,0)$ and its
  face value is $a$.
\end{remark}
The next lemma connects the face value of a central triangle with entries of
normal vectors:

\begin{lemma}\label{MMM}
  Let $\bigtriangleup$ be an empty central triangular face with three adjacent faces $\bigtriangledown_i$.
  Write $n_i \coloneqq n_{\bigtriangleup,\bigtriangledown_i}$,
  $\normalvector{\bigtriangledown_i}=(a^{(i)}_1,a^{(i)}_2,a^{(i)}_3)$ ($1\leq i\leq 3$); and
  $\normalvector{\bigtriangleup}=(a_1,a_2,a_3)$. Then
  \begin{equation}\label{MM}
    -m_{\bigtriangleup}=   a_1a_2a_3\cdot \left(
      e_{\bigtriangleup}+ \sum_{i=1}^3 \frac{a^{(i)}_i}{n_ia_i} \right).
  \end{equation}
\end{lemma}

\begin{proof}
Let $(0,p_2,p_3)$ be a vertex of $\bigtriangleup$. Then $m_{\bigtriangleup}=p_2a_2+p_3a_3$. By
\eqref{eq:20} (and by a sign check) $a^{(2)}_1a_2-a^{(2)}_2a_1=n_2p_3$ and
$a^{(3)}_1a_3-a^{(3)}_3a_1=n_3p_2$.  Use these identities and \eqref{eq:10}.
\end{proof}

\begin{remark}\label{rem:prop}\mbox{}
  \begin{enumerate}[(a)]
  \item\label{item:44} Of course, all the results proved for $G(\Gamma_{\min}(f))$
    can be transformed into properties of $G^o=G^o(\Gamma_{\min}(f))$. E.g.,
    \ref{mingood} reads as follows.  The diagram $G^o$ reflects the shape and
    adjacency relations of $\Gamma_{\min}(f)$: the vertices of $G^o$ correspond to
    the faces of $\Gamma_{\min}(f)$.  The edges of $G^o$ connecting vertices
    correspond to edges of $\Gamma_{\min}(f)$ not lying in $\partial \Gamma_{\min}(f)$, and the
    legs of $G^o$ correspond to the primitive segments lying on the edges of
    $\partial \Gamma_{\min}(f)$.  By~\ref{negy}, all the leg-decorations are greater than
    \(1\), and they satisfy the divisibility properties of~\ref{prop:2}.

    Moreover, by the very definition, the intersection orbifold matrix $I^o$
    (cf.~\ref{gr12}) can be read from $G^o$ and also the combinatorial areas
    ${\{g(\bigtriangleup)\}}_{\bigtriangleup\in\calf_c}$, needed in~\eqref{eq:12}.  Indeed, $g(\bigtriangleup)+2 $ equals
    the degree of the corresponding vertex in $G^o$ (cf.~\eqref{eq:3} and
    \eqref{eq:rhs}).
  \item\label{item:45} One may ask: how easily can $\Gamma_{\min}(f)$ be recognized
    from $G^o$? Well, rather hardly! Already the types
    (\ref{item:27}--\ref{item:29}) of \ref{cor:str} are hardly recognizable.
    The `easy' cases include $\nrt_*$ (cf.~\ref{symbols}) recognized via
    \ref{lem:4}, or $\ctr_3$ ($G^o$ has a vertex adjacent with three other
    vertices). In these cases, one also recognizes the vertices corresponding
    to hands (vertices adjacent with one vertex), or to central faces. Also,
    if $G^o$ has one vertex, then it corresponds to a one-faced diagram. But
    all the other families cannot be easily separated. E.g., it is hard to
    separate the case of a central triangle with two non-degenerate arms from
    the case of a central edge. In these cases, it is not easy at all to find
    the hands or central triangles.

    Another difficulty arises as follows. Consider an arm (with many
    triangles) crossing, say, the \(z\sb{3}\) axis. There are two types of
    triangles in it, depending on whether the non-crossing edge is on the
    $z_1z_3$ or $z_2z_3$ plane. These types are invisible from the
    \emph{shape} of $G^o$ (and will be determined using technical arithmetical
    properties of the decorations).
  \end{enumerate}
\end{remark}

\section{Starting the inverse algorithm}\label{inv}

\subsection{The main result} \label{sec:main-results}

\begin{para}\label{M11}
  We consider analytic germs $f\colon (\setC^3,0) \to (\setC,0)$ with \emph{isolated
    singularity} at \(0\), with \emph{non-degenerate Newton principal part},
  and with \emph{rational homology sphere} link $K(f)$. At a combinatorial
  level, this means that we consider all the Newton boundaries $\Gamma$ with
  \eqref{eq:2} and \eqref{eq:rhs}. Oka's algorithm §\ref{OKAsect} provides a
  resolution graph $G(\Gamma)$ for each such $\Gamma$.  Such a graph, in general, is not
  good minimal. But \ref{mingood} guarantees that if some graph $G$ can be
  obtained by this procedure, then also the good minimal resolution graph
  $G_{\min}$ associated with $G$ (obtained from $G$ by repeated blow downs of
  $(-1)$-vertices of degree less than \(3\)) can be obtained by running the
  algorithm for a minimal representative $\Gamma_{\min}$ of $\Gamma$.

  Recall that all the resolution graphs which are equivalent modulo blowing
  up/down $(-1)$-vertices can be regarded as the plumbing graphs of the same
  plumbed \(3\)-manifold, the link $K(f)$. By \cite{NeuCalc}, this class of
  graphs, and the unique good minimal one, can be recovered from the oriented
  topological type of the link $K(f)$.

  Recall also that to any graph $G$ one can associate the orbifold diagram
  $G^o$.

  Our next result, which also implies Theorem \ref{th:1} from the
  introduction, says that Oka's algorithm is, basically, injective:
\end{para}
\begin{theorem}\label{th:2}
  The d-minimal representative $\Gamma_{\min}$ (up to a permutation of coordinates)
  can be uniquely recovered from the orbifold diagram $G^o$ associated with
  $G(\Gamma_{\min})$.
\end{theorem}

\begin{corollary}\label{okaoka}
  We consider germs as in \ref{M11}. Then one has:
  \begin{enumerate}[(i)]
  \item\label{item:46} The orbifold diagram associated with the good minimal
    resolution contains the same information as $G_{\min}$ itself.
  \item\label{item:47} If the links of two germs $f_0$ and $f_1$ are
    homeomorphic then there exist germs ${\{g_i\}}_{i=0}^k$ (as in \ref{M11})
    and a coordinate-permutation $\sigma$ so that $g_0=f_0$, $g_k=f_1\circ \sigma$, and
    $g_i+t(g_{i+1}-g_i)$ ($0\leq i<k$) is a $\mu$-constant deformation
    corresponding to one of the moves $\jpn$.
  \end{enumerate}
\end{corollary}

\begin{para}[Outline of the algorithm]\label{alg:main-part}  The inverse algorithm
  which recovers the d-minimal representative from $G^o$ is rather long.  It
  distinguishes \(3\) cases depending on the number \(N\) of nodes:

  \begin{description}
  \item[\hypertarget{N0}{N=0}] The minimal resolution graph has \emph{no
      nodes}.  This is the simplest case solved in \ref{cor:3}: Let $n$ denote
    the determinant (equivalently, $G^o$ is a `free' edge with decoration
    $n$).  The representative is the diagram of $z_1^n+z_2z_3$.
  \item[\hypertarget{N1}{N=1}] The minimal resolution graph is
    \emph{star-shaped} (equivalently, $G^o$ has only one vertex).  This case
    corresponds to (deformations of) isolated weighted homogeneous germs.
    Then \cite{MN} proves that from the resolution graph one can recover the
    supporting plane $\pi$ of the unique face of the (representative) Newton
    boundary (or the weights). Next we provide an even shorter argument. By
    \cite{Pi}, the Poincaré series of the graded algebra of the germ $f$ is
    recovered from $G(f)$. But this is a rational function of type
    $(t^m-1)/((t^{a_1}-1)(t^{a_2}-1)(t^{a_3}-1))$ codifying the equation $\sum
    a_iz_i=m$ of $\pi$, cf.~\cite{Wa}. Putting all the possible lattice points
    on $\pi$, we get the canonical representative of $\Gamma$.

    A long combinatorial case by case verification recovers $\pi$ from $G^o$,
    too, which the patient reader may rediscover using the classification in
    Appendix~\ref{WHC}.
  \item[\hypertarget{Ngreater1}{N\textgreater{}1}] The orbifold diagram has
    \emph{at least two vertices}.  This is the subject of the remaining
    sections, an outline of it is given here.  The procedure involves three
    main technical steps:
    \begin{enumerate}[(1)]
    \item\label{item:48} \emph{arm preprocessing} which provides partial
      information about the arms and about the face(s) behind the shoulders
      (§\ref{sec:deal-with-arm});
    \item\label{item:49} determination of the \emph{center} (§§
      \ref{sec:triangle-centre}, \ref{sec:ianrt} and \ref{remaining});
    \item\label{item:50} \emph{arm postprocessing} which calculates the arms
      completely (§\ref{sec:arm-postprocessing}).
    \end{enumerate}
    In fact, (\ref{item:48}) for an arm runs only if we \emph{know the
      position of the corresponding hand}, otherwise it should be preceded by
    a \emph{hand-search} step.
  \end{description}
  In the next two subsections, we are going to discuss the easier arm
  preprocessing and postprocessing, respectively.  They are uniform no matter
  how the diagram looks like.

  On the other hand, we devote more sections to determine the center, since it
  distinguishes many cases depending on how the center looks like and uses a
  separate algorithm in every case.
\end{para}

\subsection{Arm preprocessing}\label{sec:deal-with-arm}

\begin{definition} \label{def:3} Let us consider a non-degenerate arm of a
  d-minimal Newton diagram in the direction, say, of the \(z\sb{3}\) axis. Its
  \emph{basic data} consists of the following:
  \begin{enumerate}[(1)]
  \item\label{item:51} the correspondence $\kappa$ between the triangles and edges
    of the arm and the corresponding vertices and decorated edge groups of
    $G^o$, respectively;
  \item\label{item:52} the first and second coordinates of the vertices of the
    triangles of the arm;
  \item\label{item:53} the third coordinates of the normal vectors of all the
    triangles of the arm and also of the (compact or non-compact) face of
    $\Gamma_+$ opposite the shoulder;
  \item\label{item:54} the face values of the non-compact faces adjacent to
    the triangles of the arm associating these numbers to the corresponding
    `half-free' edges of $G^o$.
  \end{enumerate}
\end{definition}
The basic data is an invariant of the arm, i.e.\ it is independent of the
parts of $\Gamma$ outside the arm, and also does not depend on the choice of
coordinates; explicit coordinates are used in the definition only for
simplicity of language.  In particular, in the language
of~(\ref{item:51})--(\ref{item:54}) above, it is only well-defined up to a
permutation of the first two coordinates. Nevertheless, this permutation is
\emph{global}: if we exchange the coordinates in one triangle, then we have to
exchange in all of them.  (In~(\ref{item:54}) the face values are independent
of the permutation of the first two coordinates. In fact, they are \(0\)
excepting maybe one leg of the hand.)

It is convenient to distribute the basic data among the triangles: The
\emph{basic data} of a triangle $\bigtriangleup$ of an arm (in the direction of the
\(z\sb{3}\) axis) consists of the first and second coordinates of its vertices
and the third coordinate of its normal vector, and also the correspondence $\kappa$
between the edges of $\bigtriangleup$ in $\Gamma$ and edge groups of $\kappa(\bigtriangleup)$ in $G^o$. The basic
data of a triangle is part of the basic data of the containing arm; the choice
of coordinates agree with the choice of coordinates for the arm. (If a
triangle is contained in several arms then it has a separate basic data for
each of the containing arms.)

\begin{para}[The aim of arm preprocessing] \emph{ Assume that we identify in
    $G^o$ the vertex corresponding to the hand of an arm in the d-minimal
    representative. The aim of arm preprocessing is to determine from $G^o$
    the basic data of this arm}. We will compute the basic data of the
  triangles of the arm one after the other beginning at the hand. Meanwhile,
  we will also recognize when we reach the shoulder of the arm, and we will
  compute the third coordinate of the normal vector of the next face as
  required by~\ref{def:3}(\ref{item:53}).
\end{para}
\begin{para}
  We start with the basic data of the hand $\bigtriangleup$. At this stage we are free to
  make any choice of coordinates: we assume that the arm is in the direction
  of \(z\sb{3}\) axis; and if $\bigtriangleup$ has any edge with interior lattice points
  (say, $t-1$ of them) then this edge sits on the \(z\sb{2} z\sb{3}\) plane.
  Let $\kappa(\Delta)$ be the corresponding vertex in $G^o$. The next paragraph collects
  some facts about \emph{decorations of the legs} adjacent to $\kappa(\Delta)$.
  By~\ref{negy}, all of them are greater than \(1\).
\end{para}

\begin{para}\label{facts-hands} If $\bigtriangleup$ intersects the
  \(z\sb{3}\) axis, then there are two types of leg-decorations
  (cf.~\ref{mingood}): $t$ legs decorated with
  $\innerdeterminant{\bigtriangleup}{\coordinatevector{1}}$ and one leg with
  $\innerdeterminant{\bigtriangleup}{\coordinatevector{2}}$. Notice that by
  \ref{prop:2}(\ref{item:39}) one has
  $\mygcd{\innerdeterminant{\bigtriangleup}{\coordinatevector{1}},
    \innerdeterminant{\bigtriangleup}{\coordinatevector{2}}}=1$, and by \eqref{eq:p2} the
  vertices of $\bigtriangleup$ have the form $(0,0,*)$, $(0,
  t\innerdeterminant{\bigtriangleup}{\coordinatevector{2}},*)$, and
  $(\innerdeterminant{\bigtriangleup}{\coordinatevector{1}},0,*)$.  This last fact together
  with \eqref{eq:22} implies $a_{3}=
  \innerdeterminant{\bigtriangleup}{\coordinatevector{1}}\cdot
  \innerdeterminant{\bigtriangleup}{\coordinatevector{2}}$.

  Otherwise, if $\bigtriangleup$ has an edge of type \([(a,0,c),(0,1,b)]\) ($a>0$), set
  \(\namedvector{n} \coloneqq \vectorbycoordinates{1,a,0}\) as in
  \ref{lem:10}. Then there are $t$ legs with decoration
  $\innerdeterminant{\bigtriangleup}{\coordinatevector{1}}$ and one leg with
  $\innerdeterminant{\bigtriangleup}{\namedvector{n}}$. By \eqref{eq:23} we have
  $\innerdeterminant{\bigtriangleup}{\namedvector{n}}= a_{3}$, by \eqref{eq:p2} we have
  $\innerdeterminant{\bigtriangleup}{\coordinatevector{1}}=a$, and by \eqref{eq:24} we have
  $\innerdeterminant{\bigtriangleup}{\coordinatevector{1}} \mid
  \innerdeterminant{\bigtriangleup}{\namedvector{n}}$ (since $r=0$).  Again by
  \eqref{eq:24}, the vertices of $\bigtriangleup$ have the form $(0,1,*)$,
  $(\innerdeterminant{\bigtriangleup}{\coordinatevector{1}},0,*)$ and $(0,1+ t
  \innerdeterminant{\bigtriangleup}{\namedvector{n}}/
  \innerdeterminant{\bigtriangleup}{\coordinatevector{1}},*)$.  Notice that
  $\innerdeterminant{\bigtriangleup}{\namedvector{n}}=
  \innerdeterminant{\bigtriangleup}{\coordinatevector{1}}$ may happen only in the case
  \ref{prop:2}(\ref{item:38}).

  For the corresponding face values, see \ref{lem:11b}.
\end{para}

\begin{para}[Algorithm: the basic data of a hand]\label{alg-hand}
  Let $N$ be the set of all the decorations of the legs adjacent to $\kappa(\Delta)$.
  One may have the following situations:
  \begin{enumerate}[(a)]
  \item\label{item:55} $N=\{n_1,n_2\}$ with $\mygcd{n_1,n_2}=1$. One of them,
    say $n_2$, decorates exactly one leg; the other one, $n_1$, decorates
    several legs, say $t$ of them. (If $t=1$ then the construction is
    symmetric.) Then (up to a permutation of the coordinates, cf.~\ref{def:3})
    the vertices of the hand have the form $(0,0,*)$, $(0,tn_2,*)$ and
    $(n_1,0,*)$; and the third coordinate of $\normalvector{\bigtriangleup}$ is $n_1n_2$.
    All the face values are \(0\).
  \item\label{item:56} $N=\{n_1,n_2\}$ with $n_1 \mid n_2$. Then the number of legs
    decorated by $n_1$ will be denoted by $t$, and (automatically) $n_2$
    decorates one leg. The face value of legs with $n_1$-decoration is \(0\),
    but the face value of the unique $n_2$-decorated leg is $n_1$. The
    vertices of the hand have the form $(0,1,*)$, $(n_1,0,*)$ and $(0,1+ t
    n_2/ n_1,*)$ and the third coordinate of $\normalvector{\bigtriangleup}$ is $n_2$.
  \item\label{item:57} $N=\{n\}$. Then let $t+1$ be the total number of legs. We
    set $n_1=n_2=n$. Then the basic data of the hand is given by the same
    formulas as in (\ref{item:56}). We separate one leg (with face value $n$),
    the others form another group (with face value \(0\)).
  \end{enumerate}
\end{para}

\begin{para}[Arm continuation]\label{sec:arm-continuation}
  Assume that we have computed from $G^o$ the basic data of the triangles
  $\bigtriangleup_1,\dotsc, \bigtriangleup_k$ (belonging to an arm in the direction of $z_3$, where
  $\bigtriangleup_1$ is the hand, and $\bigtriangleup_i$ is adjacent to $\bigtriangleup_{i+1}$) in such a way that
  the coordinate-ambiguities are compatible (i.e., if we fixed coordinates for
  $\bigtriangleup_1$, then for all the other $\bigtriangleup_i$ we respect the same choice).  Our aim is
  \emph{to determine the part of the basic data corresponding to the next
    face}, i.e.\ the third coordinate of its normal vector, and whether it
  belongs to the arm.  If yes then we also compute its basic data.

  We write \(\bigtriangleup\coloneqq \bigtriangleup_k\) and set $\kappa(\bigtriangleup)$ for the corresponding vertex of
  $G^o$. By the inductive step, we have already computed the correspondence
  $\kappa$ of all the edges of \(\bigtriangleup\) with the edges adjacent to $\kappa(\bigtriangleup)$. Let \(\bigtriangledown\)
  be the next face of $\Gamma_+$, adjacent to \(\bigtriangleup\), and set $\gamma\coloneqq \bigtriangledown\cap \bigtriangleup$. The
  face \(\bigtriangledown\) is compact if and only if $\kappa(\gamma)$ connects two vertices of $G^o$,
  one of them is obviously $\kappa(\bigtriangleup)$. If this is the case, we set \(\kappa(\bigtriangledown)\) for
  the other end.

  In any situation, we need the third coordinate $a_3$ of the normal vector
  $\normalvector{\bigtriangledown}$. This can be computed from $G^o$ and the basic data of
  the triangles $\bigtriangleup_i$ using~\eqref{eq:10}.

  Next, if \(\bigtriangledown \) is compact, we wish to decide whether it belongs to the arm.
  The face \(\bigtriangledown \) is a non-removable trapezoid if and only if $\kappa(\bigtriangledown)$ admits
  four distinguishable groups of adjacent edges (cf.~\ref{lem:4}). In this
  case, clearly, \(\bigtriangledown \) does not belong to the arm. The same is true if $\kappa(\bigtriangledown)$
  has no legs (which happens if and only if $\bigtriangledown$ is a central triangle).
  Therefore, assume that $\bigtriangledown$ is a triangle with at least one adjacent leg.

\begin{lemma}\label{oszthatosag} The face $\bigtriangledown$
  belongs to the arm if and only if $\kappa(\bigtriangledown)$ has an adjacent leg whose
  decoration divides $a_3$.
\end{lemma}
\begin{proof} Write $\gamma =AB$, and let $C$ be the third vertex of $\bigtriangledown$. If $C$ is
on the $z_2z_3$ or $z_1z_3$ planes, then the leg associated with $BC$ or $AC$
divides $a_3$ by~\ref{prop:2}(\ref{item:35}).  If $C$ is not sitting on one of
these two planes, then it can have a leg only if at least one of $AC$ and $BC$
is a crossing edge. But, by~\ref{prop:2}(\ref{item:36}), such an edge
determinant divides $a_3$ only if $\bigtriangledown$ is a moving triangle whose moving vertex
$C$ can be moved to a coordinate axis. But this would contradict the
definition of the d-minimal representatives in~\ref{cor:str2}.
\end{proof}

If $\bigtriangledown$ does not belong to the arm, then we stop (having all the basic data of
the arm).

Next, assume that $\bigtriangledown$ belongs to the arm. Then we have to identify edges of
$\bigtriangledown$ with edges adjacent to $\kappa(\bigtriangledown)$ in $G^o$, and to determine the first two
coordinates of $C$.
 
First we identify the leg-decoration $\legdeterminant{\bigtriangledown}$, adjacent to $\kappa(\bigtriangledown)$,
which corresponds to the edge $\alpha$ of $\bigtriangledown$ which lies on a coordinate plane.  By
\ref{prop:2}(\ref{item:35}), it divides $a_3$. We claim that
\emph{$\legdeterminant{\bigtriangledown}$ is the largest leg-decoration adjacent to $\kappa(\bigtriangledown)$
  which divides $a_3$}. Indeed, we have to check only the case when $\kappa(\bigtriangledown)$ has
two leg-decorations $N=\{n_1,n_2\}$ (the determinants of $AC$ and $BC$, one edge
sitting on a coordinate plane, the other being a crossing edge), both dividing
$a_3$. Then, by~\ref{prop:2}, $\bigtriangledown$ is a moving triangle such that $C$ can be
moved to both coordinate axes, and by the construction of the d-minimal
representative (cf.~\ref{cor:str2}), the determinant of the edge which lies on
the coordinate plane is the larger one.

Now, we fix an edge of $\bigtriangleup$ (whose determinant will be denoted by
$\legdeterminant{\bigtriangleup}$) which lies on a coordinate plane (which is either
$z_1z_3$ or $z_2z_3$ determined clearly by the basic data of $\bigtriangleup$). Denote this
plane by $\pi$.  Then, by~\ref{lem:11}, $\alpha$ lies on $\pi$ if and only if
$\legdeterminant{\bigtriangleup} = \legdeterminant{\bigtriangledown} \mid \innerdeterminant{\bigtriangleup}{\bigtriangledown}$.

This is valid for $C$, too, hence this clarifies whether $C$ is on the
$z_1z_3$ or $z_2z_3$ plane. Finally, we have to compute the first two
coordinates of $C$. One of them is \(0\) (depending whether $\pi$ is the
$z_1z_3$ or $z_2z_3$ plane), the other can be determined using~\eqref{eq:22}.

Then we add $\bigtriangledown$ to the triangles ${\{\bigtriangleup_i\}}_i$ and repeat the `arm continuation'
process by induction.
\end{para}

\subsection{Arm postprocessing}
\label{sec:arm-postprocessing}

\begin{para}
  The \emph{arm postprocessing} step assumes the knowledge of two sets of
  data: the first one is the basic data coming from arm preprocessing, the
  second one is some knowledge about the face on the other side of the
  shoulder (which is usually the center but not always). More precisely, let
  $\bigtriangledown$ be the face in $\Gamma_+$ containing the shoulder of the arm, but not
  contained in the arm. (It is non-compact if and only if the edge $\kappa(\bigtriangledown)$ in
  $G^o$ is a leg.)

  \begin{enumerate}[I.]
  \item\label{item:58} \itemtitle{The first set of data:} the basic data of an
    arm in the direction of the \(z\sb{3}\) axis modulo the ambiguity of a
    permutation of the first two coordinates.  Recall that the basic data of
    the arm also includes the knowledge of
    $\scalarproduct{\normalvector{\bigtriangledown}}{\coordinatevector{3}}$.
  \item\label{item:59} \itemtitle{The second set of data:} consists of all the
    coordinates of $\normalvector{\bigtriangledown}$ and of the shoulder of the arm in some
    choice of coordinates \(z_1\) and \(z_2\). The coordinate \(z_3\) is the
    same as in the first set.
  \end{enumerate}
  Hence, by assumption, we have a `half-compatibility' connecting the two
  choices of coordinates in the two sets of data: the third coordinate $z_3$
  in the basic data (\ref{item:58}) and for the pair
  $(\normalvector{\bigtriangledown},\text{shoulder})$ in (\ref{item:59}) are matched. But, a
  priori, we do not know how to identify the other (i.e.\ the first two)
  coordinates in the two sets of data.

  The \emph{aim of arm postprocessing} is (using (\ref{item:58}),
  (\ref{item:59}) and $G^o$) to compute all the coordinates of the vertices of
  the arm in a unified choice of coordinates for the sets of data.
\end{para}

\begin{para}[Unifying the first two coordinates]\label{match}
  Let \(z_1\), \(z_2\), \(z_3\) be the coordinates in which we describe the
  second set of data (\ref{item:59}). Assume that in these coordinates the end
  points of the shoulder are $A=(0,p_2,p_3)$ and $B=(q_1,0,q_3)$. (Notice that
  the data (\ref{item:58}) recognizes the first two coordinates up to their
  permutation of $A$ and $B$. Hence, if $p_2\neq q_1$ then these information
  already unifies the coordinates. But, in general, we have to do more.)

  Let $\bigtriangleup$ be the last triangle of the arm (i.e.\ $\bigtriangleup\cap\bigtriangledown=AB$). Fix an edge $\alpha$ of
  $\Delta$ in $\Gamma$ whose edge-group $\kappa(\alpha)$ in $G^o$ contains, say, $t$ legs. We will
  determine whether $\alpha$ is on the $z_1z_3$ or $z_2z_3$ plane: this will orient
  all the basic data (\ref{item:58}) in accordance with ${\{z_i\}}_i$.

  First we compute (in the coordinates ${\{z_i\}}_i$) the normal vector
  $\normalvector{\bigtriangleup}$. For this, notice that $\innerdeterminant{\bigtriangleup}{\bigtriangledown} \cdot
  \overrightarrow{AB} = \normalvector{\bigtriangleup} × \normalvector{\bigtriangledown}$.  This follows
  from~\ref{lem:9}(\ref{item:77}) up to a sign; the sign is a consequence of
  the right-hand rule for vector products.  Since $\overrightarrow{AB}$,
  $\normalvector{\bigtriangledown}$ and $\innerdeterminant{\bigtriangleup}{\bigtriangledown}$ are known, this identifies
  $\normalvector{\bigtriangleup}$ up to a summand of a scalar multiple of
  $\normalvector{\bigtriangledown}$. Since
  $\scalarproduct{\normalvector{\bigtriangleup}}{\coordinatevector{3}}$ is also known (from
  (\ref{item:58})) and
  $\scalarproduct{\normalvector{\bigtriangledown}}{\coordinatevector{3}}$ is not \(0\), these
  facts determine $\normalvector{\bigtriangleup}$ completely.

  Now, we determine whether $\alpha$ lies on $z_1z_3$ or $z_2z_3$ plane. Recall
  that from the basic data (\ref{item:58}), we know the set of the first two
  coordinates of $\alpha$: one of them is \(0\), the other one is, say, $\delta(\alpha)>0$.
  E.g., if $\alpha=[(q_1,0,q_3),(q_1',0,q_3')]$, then $\delta(\alpha)=q_1-q_1'$, and $\delta(\alpha)/t$
  is a positive \emph{integer}, known from the basic data (\ref{item:58}).
  Then \eqref{eq:22} and \eqref{eq:ijk}, for $i\in\{1,2\}$, reads as:
  \begin{equation*}
    \alpha\subset\text{$z_iz_3$ plane} \implies
    a_{3}= \mygcd{ a_{3}, a_{i}}\cdot\delta(\alpha)/t.
  \end{equation*}
  Since $\normalvector{\bigtriangleup}$ is primitive, and the above greatest common divisor
  is, in fact, a leg-determinant, hence it is greater than \(1\)
  by~\ref{negy}, the right-hand side of the above identity cannot be true for
  both $i=1,2$ simultaneously. This fact determines which coordinate plane
  contains $\alpha$.
\end{para}

\begin{para}[The complete determination of the arm]
  Now, using \ref{match}, we can write all the basic data (\ref{item:58}) in
  the coordinates \(z_1\), \(z_2\), \(z_3\) of (\ref{item:59}). Notice that
  the basic data (\ref{item:58}) determines completely all the normal vectors
  and all the face values associated with the non-compact faces adjacent to
  the arm (cf.~\ref{lem:11b}). Moreover, from (\ref{item:59}) we know the
  normal vector and the face value of $\bigtriangledown$. Hence the affine equations of all
  the triangles in the arm follow from the systems \ref{cor:negdef} (where
  ${\overline{\calf}}_c$ is the index set of triangles of the arm).
\end{para}

\subsection{The complete inverse algorithm for $\ctr_3$.}
\label{sec:triangle-centre}

We end this section by clarification of case $\ctr_3$. First notice that this
family can be identified using the diagram $G^o$: it has a (unique) vertex $v$
with three adjacent vertices and without any legs. In this subsection we
assume that $G^o$ has this property.

The vertex $v$ corresponds to the central triangle $\bigtriangleup$. The other vertices can
be grouped in three, each group consisting of a string of adjacent vertices
corresponding to the three arms of the diagram. The hands correspond to
vertices with exactly one adjacent vertex.

We mark the three vertices corresponding to the hands (or directions of the
arms) with the three coordinates. Here we are free to make any marking (up to
a permutation of the coordinates).  We fix one. Once this choice is made, let
us denote the coordinates of $\bigtriangleup$ by \((0,p_2,p_3)\), \((q_1,0,q_3)\) and
\((r_1,r_2,0)\). At this stage these entries are unknowns.

Now, we preprocess the arms.  E.g., for the arm in the direction of $z_3$ we
obtain the basic data of that arm up to a permutation of the first two
coordinates. In particular, we obtain
\begin{enumerate}[(a)]
\item\label{item:60} the first two coordinates of the shoulder
  $[(0,p_2,p_3),(q_1,0,q_3)]$ up to a permutation, hence the set
  $\calS=\{q_1,p_2\}$;
\item\label{item:61} $\scalarproduct{\normalvector{\bigtriangleup}}{\coordinatevector{3}}$;
\item\label{item:62} the face values of the legs of this arm.
\end{enumerate}
Summing up for all three arms, we get
\begin{enumerate}[(A)]
\item\label{item:63} the pairs of coordinates \(\{ r_2, q_3 \}\), \(\{ q_1, p_2
  \}\), \(\{ p_3, r_1 \}\);
\item\label{item:64} $\normalvector{\bigtriangleup}$;
\item\label{item:65} the face values of all the legs and, hence, the face
  value $\facevalue{\bigtriangleup}$ of $\bigtriangleup$, too, from the system \eqref{eq:12}
  via~\ref{cor:negdef}.
\end{enumerate}
\begin{lemma}\label{UJUJ}
  Let \((0,p_2,p_3)\), \((q_1,0,q_3)\) and \((r_1,r_2,0)\) be the vertices of
  an empty triangle $\bigtriangleup$ (with $p_2, p_3, q_1, q_3, r_1, r_2>0$). Then these
  coordinates are uniquely determined by:
  \begin{enumerate}[(i)]
  \item\label{item:66} the sets \(\{ r_2, q_3 \}\), \(\{ q_1, p_2 \}\), \(\{ p_3,
    r_1 \}\);
  \item\label{item:67} the normal vector $(a_1,a_2,a_3)=\normalvector{\bigtriangleup}$ of
    $\bigtriangleup$;
  \item\label{item:68} the face value $\facevalue{\bigtriangleup}$ of $\bigtriangleup$.
  \end{enumerate}
\end{lemma}
Coming back to our original situation, \ref{UJUJ} determines the central
triangle $\bigtriangleup$. Then we postprocess the arms to calculate all the missing data
about $\Gamma$.

\begin{proof}[Proof of~\ref{UJUJ}.]
We have a $\setZ_{2}$-ambiguity for each set of~(\ref{item:66}) and we wish to
select the correct choice from the \(2^3\) possibilities.  For this, first
assume that we are able to decide which element of the set \(\{ r_2, q_3 \}\) is
\(r_2\) and which one is \(q_3\).  Then we claim that the other two
ambiguities disappear. Indeed, the face value identities written for the
vertices of $\bigtriangleup$
\begin{equation}\label{eq:face10}
  a_1q_1+a_3q_3=a_2p_2+a_3p_3=a_1r_1+a_2r_2=\facevalue{\bigtriangleup},
\end{equation}
(where $a_k>0$ for all $k$) and (\ref{item:66}) provide $\bigtriangleup$. Hence, we have at
most two choices: either the correct one for every arm or the wrong one for
every arm (i.e.\ when we interchange $r_2$ with $q_3$ and $q_1$ with $p_2$ and
$p_3$ with $r_1$). We claim that the wrong choice can be ruled out.  Indeed,
assume that both choices of system of integers satisfy the
formula~\eqref{eq:19} for $\normalvector{\bigtriangleup}$ and \eqref{eq:face10}. Notice
that without loss of generality, we may assume that \(r_1\) is the smallest
among \(r_1\), \(p_2\) and \(q_3\). Write (a part of) \eqref{eq:19} and
\eqref{eq:face10} for both choices:
\begin{align} \label{eq:14}
  p_2 q_3 + r_2 p_3 - r_2 q_3 &= q_1 r_2 + q_3 r_1 - q_3 r_2 = a_1,\\
  \label{eq:15}
  q_3 r_1 + p_3 q_1 - p_3 r_1 &= r_2 p_3 + r_1 p_2 - r_1 p_3 = a_2,\\
  \label{eq:face20}
  a_2p_2+a_3p_3=a_1r_1+a_2r_2&=\facevalue{\bigtriangleup}=a_2q_1+a_3r_1=a_1p_3+a_2q_3.
\end{align}
Then \eqref{eq:14} implies that $ 0 \leq (p_2 - r_1)/r_2 =(q_1 - p_3)/q_3$, hence
\(q_1 \geq p_3\) too. Thus:
\begin{equation*}
  \label{eq:16}
  a_2 = q_3 r_1 + p_3 q_1 - p_3 r_1 \geq r_1^2 + p_3^2 - r_1 p_3 = (r_1 - p_3)^2
  + r_1 p_3 > \lvert r_1 -p_3 \rvert.
\end{equation*}
On the other hand, from \eqref{eq:face20} expressing \(a_1\) and \(a_3\)
yields: \(a_1 = a_2(q_3 - r_2)/(r_1 - p_3)\) and \(a_3 = a_2(p_2 - q_1)/(r_1 -
p_3)\).  But this contradicts to the fact that \(\normalvector{\bigtriangleup}\) is
primitive:
\begin{equation*}
  \mygcd{a_1,a_2,a_3} = \frac{a_2}{\lvert r_1 - p_3 \rvert} \mygcd{q_3 - r_2, r_1
    - p_3, p_2 - q_1} > 1.
\end{equation*}
\end{proof}

\section{The inverse algorithm for the families $\nrt_*$}
\label{sec:ianrt}

\subsection{The start}

\begin{para} By~\ref{lem:4}, the family $\nrt$ can be identified from $G^o$:
  it has a (unique) vertex $v$ with four different types of edges.  In this
  section we assume that $G^o$ has this property.

  The vertex $v$ corresponds to a non-removable trapezoid $\bigtriangleup$. This vertex
  always has at least one leg group (corresponding to the bottom edge).
\end{para}

\begin{lemma} \label{lem:7} The diagram has at least one non-degenerate arm.
\end{lemma}
\begin{proof} Assume that the top edge is the shoulder of a degenerate arm.
Then, write the coordinates of the vertices $\bigtriangleup$ as in \ref{lem:1}. Then (up to
a permutation of the first two coordinates) $q=1$, and by
\ref{ex:quadrangle-ambiguity}, $r_2>0$ and $n>1$. For $r_2=1$, the trapezoid
$\bigtriangleup$ can be enlarged and is removable. Hence $\Gamma$ has a non-degenerate arm in
the direction of $z_1$.
\end{proof}

The hands can also be identified in $G^o$: they are those vertices (different
from $v$) which have one adjacent vertex.  The next algorithm splits according
to their number.

Notice also, that by \ref{lem:8}, and with the notation of \ref{lem:1}, the
normal vector of $\bigtriangleup$ is
\begin{equation}
  \label{eq:17}
  \normalvector{\bigtriangleup} = \vectorbycoordinates{nq,np, r_1 q +  r_2 p + (t - 1) pq}.
\end{equation}

\subsection{The case of three non-degenerate arms: $\nrt_3$.}
\label{sec:quadr-cent-one}

This case has many similarities with \ref{sec:triangle-centre}; but, in fact,
it is simpler since the legs of $v$ help in the procedure. Let the decoration
of the unique leg group of $v$ be $d$.

We start by preprocessing the three arms. This provides the coordinates of
$\normalvector{\bigtriangleup}$ (up to a permutation). Since we already identified the
vertex of the central face, we know when we arrive to the shoulder.
Nevertheless, at this step, we see a difference between the side-arms and the
top-arm.  Consider e.g.\ a side arm and the `hidden' triangle (as part of $\bigtriangleup$)
formed by the shoulder and the base edge. With this triangle the \emph{arm
  continuation procedure \ref{sec:arm-continuation} is not obstructed}, in
other words, (by~\ref{oszthatosag}) $d$ divides the corresponding coordinate
of $\normalvector{\bigtriangleup}$.  For the top-arm this is not the case.

Therefore, $d$ divides exactly two coordinates of $\normalvector{\bigtriangleup}$. We
attach the coordinate $z_3$ to the arm for which this divisibility does not
hold (in this way its shoulder will be the top edge and the bottom edge will
sit on the \(z\sb{1} z\sb{2}\) plane).  The coordinates $z_1$ and $z_2$
(chosen arbitrarily) will be attached to the other two strings of vertices.

Since the face values of the legs of $v$ are \(0\), and all the other face
values associated with legs have been determined during arm preprocessing,
\eqref{eq:12} and \ref{cor:negdef} provide the face value of $\bigtriangleup$. In
particular, we get the equation of the affine plane supporting $\bigtriangleup$. Since the
top edge is primitive (and parallel to the \(z\sb{1} z\sb{2}\) plane), this is
enough for its identification. In particular, with the notations of
\ref{lem:1}, we get \(n\), \(p\) and $q$ (in fact, $n=d$ by \eqref{eq:17}).
During preprocessing of the arm in the direction $z_1$, we have obtained the
set $\{n,r_2\}$, but $n$ is already identified, hence we obtain $r_2$, too.
Similarly, we get $r_1$.  Thus we know all the vertices of $\bigtriangleup$, hence the
algorithm finishes by postprocessing the arms.

\subsection{The case of two non-degenerate arms: $\nrt_2$.}
\label{sec:quadr-cent-two}

\begin{para} We have two different leg groups (one of them attached to the
  bottom), and two non-degenerate arms. First we have to determine whether the
  shoulders of the non-degenerate arms are the side edges, or one of them is
  the top edge. This can be decided by the following divisibility property.
  Its proof uses \eqref{eq:17} and the fact that $\normalvector{\bigtriangleup}$ is
  primitive (the details are left to the reader):
\end{para}
\begin{lemma}\label{decide}
  Consider a trapezoid $\bigtriangleup$ with coordinates as in \ref{lem:1}. Assume that it
  has two non-degenerate arms. Then (up to permutation of the first two
  coordinates) there are two possibilities:

  \begin{enumerate}[{Case} 1.]
  \item\label{item:69} $q=1$, $n>1$, $r_1>1$ and $r_2>1$, the direction of the
    non-degenerate arms are the $z_1$ and $z_2$ axes, and
    $\normalvector{\bigtriangleup}=(a_1,a_2,a_3)=(n,np,*)$, hence $a_1 \mid a_2$.
  \item\label{item:70} $r_1=0$ or \(r_1=1\), and $n>1$, $p>1$, $q>1$, the
    direction of the non-degenerate arms are the $z_1$ and the $z_3$ axes, but
    the coordinates $(a_1,a_2,a_3)$ of $\normalvector{\bigtriangleup}$ do not satisfy any
    divisibility relation: $a_1\nmid a_3$, $a_3\nmid a_1$.
  \end{enumerate}
\end{lemma}

Thus the algorithm starts by preprocessing the two non-degenerate arms to get
two coordinates of $\normalvector{\bigtriangleup}$. If one of them divides the other then
we are in Case~\ref{item:69} above, otherwise we are in Case~\ref{item:70}.
Next we treat each case independently.

\begin{para}[Case~\ref{item:69}. Side edges as shoulders of non-degenerate
  arms]\label{sec:leg-attached-top}

  Preprocessing the arms has provided (say) the first two coordinates of the
  normal vector \(\normalvector{\bigtriangleup}\).  We name the coordinate axes so that the
  smallest of the first two coordinates of \(\normalvector{\bigtriangleup}\) is the first
  coordinate.  Then we compute \(n\) as the first coordinate and \(p\) as the
  fraction of the first two coordinates. Since $q=1$, at this point, we know
  all coordinates of the top edge. Then we end this case by the same argument
  as in \ref{sec:quadr-cent-one}.  Preprocessing of the arm in the direction
  $z_1$ has provided the set $\{n,r_2\}$.  Since $n$ is already identified, we
  obtain $r_2$. Similarly, we get $r_1$, too.  Finally, $t+1$ is the number of
  legs of $\bigtriangleup$.  Knowing all the vertices of $\bigtriangleup$, we finish by postprocessing
  the arms.
\end{para}

\begin{para}[Case~\ref{item:70}. Top edge as shoulder of a non-degenerate
  arm]\label{sec:leg-attached-side}

  Here, one may proceed in the spirit of the other cases
  \ref{sec:quadr-cent-one} and \ref{sec:leg-attached-top}, but one may use the
  following observation as well.  We may think about this situation as the
  degeneration of $\ctr_3$ (cf.~\ref{sec:triangle-centre}).  Indeed, consider
  the trapezoid $\bigtriangleup$ as in \ref{lem:1} and cut it into two triangles along
  \([(r_1+tp,r_2,0),(0,q,n)]\). Let \(\bigtriangleup_1\) denote the lower triangle (whose
  vertices are \((0,q,n)\), \((r_1+tp,r_2,0)\) and \((r_1,r_2+tq,0)\)), and
  let \(\bigtriangleup_2\) denote the upper triangle.  We may consider $\bigtriangleup_1$ as a `virtual'
  hand with two different leg groups, and $\bigtriangleup_2$ as a central triangle with two
  `genuine' and one `virtual' arms. The degeneration consists of the fact that
  $\bigtriangleup_1$ and $\bigtriangleup_2$ are in the same plane. Nevertheless, we can apply the same
  argument. The basic data of the `virtual hand', similarly as in
  \ref{alg-hand}, together with the basic data of the `genuine' arms provide
  all the data necessary to apply \ref{UJUJ} for the empty triangle $\bigtriangleup_2$.
  Therefore, we obtain $\bigtriangleup_2$ (up to a permutation of the coordinates).
  Postprocessing the two arms and completing $\bigtriangleup_2$ to a trapezoid (in its
  supporting plane) ends the procedure.
\end{para}

\subsection{The case of one non-degenerate arm: $\nrt_1$.}
\label{sec:quadr-cent-three}

The central vertex is attached to one non-degenerate arm and three leg groups.
We denote the set of decorations of these legs by $\calD$, which contains
three different elements (cf.~\ref{lem:4}).  Preprocessing of the arm provides
a coordinate of $\normalvector{\bigtriangleup}$, denoted by $A$, and two coordinates of the
shoulder (see e.g.~\ref{sec:triangle-centre}(\ref{item:60})) forming the set
$\calS$. The above discussions (and/or Appendix~\ref{arith}) determine these
three objects \(A\), \(\calD\), \(\calS\) in terms of the integers used in
\ref{lem:1} for the trapezoid $\bigtriangleup$. Basically, (up to a permutation of the
first two coordinates) there are two possibilities (depending on whether the
shoulder of the arm is a side or top edge). (For the coordinates
$(a_1,a_2,a_3)$ of $\normalvector{\bigtriangleup}$ see \eqref{eq:17}, which shows that
$\mygcd{a_3,n}=1$.  Moreover, if $r_1=0$ then $p\mid a_3$, and if $r_2=0$ then
$q\mid a_3$.)
\begin{enumerate}[I.]
\item\label{item:71} \(q=1\), \(r_2>1\), \(p>1\), \(n>1\), and \(r_1=0\) or
  \(r_1=1\), and the (minimal) diagram has an arm in the direction of the
  \(z\sb{1}\) axis. Then, $A=a_1=n$, $\calS=\{n,r_2\}$, and
  $\calD=\{n,a_3,n^{r_1}p\}$. If $r_1=0$ then $\mygcd{n,p}=1$.
\item\label{item:72} \(p>1\), \(q>1\), but $r_1$ and $r_2$ are \(0\) or \(1\),
  and the arm is in the direction of the \(z\sb{3}\) axis.  Then $A=a_3$,
  $\calS=\{p,q\}$, and $\calD=\{n,n^{r_1}p,n^{r_2}q\}$, where $\mygcd{p,q}=1$.
\end{enumerate}
The first case satisfies $A\in\calS\cap\calD$. If $A\in\calS\cap\calD$ happens in
case~\ref{item:72} then $r_1=0$ and $r_2=t=1$ (or vice versa), and $\bigtriangleup$ is a
parallelogram with two sides on coordinate planes. Hence, by a permutation of
the $z_1$ and $z_3$ axes one arrives to the situation~\ref{item:71}.

Analyzing the above data, one derives the next algorithm to recover $\bigtriangleup$ from
$G^o$.  It has two cases.

\begin{enumerate}[1.]
\item\label{item:73} If $A\in\calS\cap\calD$ then set $n=A$ and $r_2$ is the other
  element of $\calS$. The arm is in the direction of $z_1$. Two subcases may
  occur:
  \begin{enumerate}[(a)]
  \item\label{item:74} If there exists $D\in\calD\setminus \{n\}$ with $n\mid D$, then set
    $r_1=1$ and $p=D/n$;
  \item\label{item:75} If $\calD\setminus \{n\}=\{D_1,D_2\}$ with $D_1\mid D_2$, then set
    $r_1=0$ and $p=D_1$.
  \end{enumerate}
\item\label{item:76} If $A\notin \calS\cap\calD$, then we obtain $n$ as the smallest
  element of the set $\calD\setminus\{\text{all divisors of $A$}\}$.  The set $\calS$
  has two (relative prime) elements; we declare them (arbitrarily) $p$ and $q$
  (hence we will get $\bigtriangleup$ up to a permutation of the first two coordinates).
  $r_1$ and $r_2$ are determined by the fact that
  $\calD=\{n,n^{r_1}p,n^{r_2}q\}$. The arm is in the direction of the
  \(z\sb{3}\) axis.
\end{enumerate}
In this way, we recover $\bigtriangleup$ and the position of the arm in both cases, hence
the algorithm ends by postprocessing the arm.

\section{The inverse algorithm for the remaining cases
  \protect\(\ctr_1\protect\), \protect\(\ctr_2\protect\),
  \protect\(\ced_1\protect\), \protect\(\ced_2\protect\).}\label{remaining}

\subsection{Find a hand!}

\begin{para} Next we assume that $G^o$ has at least two vertices, each vertex
  has at most two adjacent vertices, and for each vertex the number of
  adjacent vertices and leg groups together is three.  An \emph{end vertex}
  has one adjacent vertex.  The diagram $G^0$ has two end vertices.  The
  minimal subgraph generated by vertices and edges connecting them is a
  string.

  Clearly, the Newton diagram has at least one hand. Compared with the
  previous cases, now it is much harder to recognize the vertices of $G^o$
  corresponding to hands (and/or centers). A hand always corresponds to an end
  vertex, but end vertices may also correspond to central triangles (e.g.\ the
  case of moving triangle), or to the last triangle (adjacent to the shoulder)
  of an arm (e.g.\ the diagram of $z_3^az_2+z_1^b+z_1^cz_2^d+z_2^e$).

  For any end vertex $v$, consider its two leg groups.  Let $t_i$ and $n_i$
  denote the number of legs and the decoration of the leg groups for $i=1,2$.
  Set $r(v)\coloneqq n_1t_2+n_2t_1$.

  By~\ref{prop:2}(\ref{item:41}), if the two decorations of an end vertex $v$
  are \emph{not} relative prime then \emph{$v$ is a hand}. We call such an end
  vertex an \emph{easily recognizable hand} (ER-hand for short).
\end{para}

\begin{lemma}\label{handhand}
  Assume that neither of the end vertices $v_1$ and $v_2$ of $G^o$ is an
  ER-hand. We mark one of them as follows. If $r(v_1)<r(v_2)$ then $v_1$ is
  marked. If $r(v_1)=r(v_2)$ then the one with greater orbifold Euler number
  is marked. If even their orbifold Euler numbers are equal, then $G^o$ has
  only two vertices (namely, \(v_1\) and \(v_2\)) and it has an isomorphism
  permuting these two vertices.  Then we mark arbitrarily one of the vertices.

  All in all, the marked vertex is always a hand (in the last case up to this
  isomorphism).
\end{lemma}

\begin{example}\label{SYMM}
  The symmetric case occurs if the Newton diagram has only four vertices
  \((0,0,2)\), \((p,0,1)\), \((0,q,1)\) and \((r_1,r_2,0)\) (satisfying \(r_1
  q + r_2 p > 2 p q\)).  Then $G^o$ has two vertices, each having two legs
  decorated by \(p\) and \(q\).  One can check that even the resolution graph
  is symmetric.  This is surprising since the Newton diagram is not symmetric
  at all: either face is a hand, the other is a moving (central) triangle.
  Nevertheless, the algorithm recovers the asymmetric Newton diagram from a
  symmetric orbifold diagram!  (Up to permutation of coordinates, this is the
  only possibility for the symmetric case, see the proof below.)
\end{example}

\begin{proof}[Proof of \ref{handhand}.]
Fix a non-degenerate arm in the direction of the \(z\sb{3}\) axis with hand
$v_1$. Then the sum of the first two coordinates of the crossing edges of this
arm strictly increases from the hand to the shoulder. For the first segment
(closest to $v_1$) it is $r(v_1)$, cf.~\ref{facts-hands}. Assume that $v_2$
corresponds to the triangle \(\trianglebyvertices{P}{Q}{R}\), with
$P=(0,p_2,p_3)$, $Q=(q_1,0,q_3)$ and $R=(r_1,r_2,0)$. We may assume that
$p_3>0$ and $q_3>0$ (otherwise $v_2$ is a hand and we have nothing to prove).
Thus, it is enough to show $q_1\leq\det(PQ)$ (and its analogue). Since by
\eqref{eq:24} this determinant is $a_2$ (of \(\trianglebyvertices{P}{Q}{R}\)),
and $a_2= q_3r_1+p_3(q_1-r_1)$ by \eqref{eq:19}, we need $q_1\leq
q_3r_1+p_3(q_1-r_1)$. By \eqref{eq:2} at least one of $p_3 $ and $r_1$ is
\(1\), hence the inequality follows. Moreover, $r(v_1)=r(v_2)$ if and only if
$v_1$ and $v_2$ are the only vertices of $G^o$, and $p_3=q_3=1$; hence the
Newton diagram is given by $(0,0,c),(0,p_2,1),(q_1,0,1),(r_1,r_2,0)$. For
this, using \ref{sec:splice-diagram}, we get $e_{v_1}\geq e_{v_2}$, with equality
if and only if $c=2$. For $c=2$ the graph is symmetric.
\end{proof}

\begin{para}[Start of the algorithm] We fix an end vertex $v_1$ which
  correspond to a hand.  We denote the other end vertex by $v_2$ (we may not
  know yet whether it is a hand).  The algorithm starts with preprocessing the
  arm with hand $v_1$.  Depending on the outcome, we continue by \ref{EGY},
  §\ref{KETTO} or §\ref{sec:two-non-degenerate}.
\end{para}

\begin{para}[The arm contains all vertices]\label{EGY}
  We assume that the arm of $v_1$ contains all the vertices of $G^o$. We fix
  the coordinates in such a way that the arm is in the direction of $z_3$.
  Then the shoulder has the form $[(r_1,0,0),(0,p_2,p_3)]$ with $p_3=0$ or
  \(p_3=1\), and $r_1>1$. Let $\namedvector{a}$ be the normal vector (of the
  non-compact face) beyond the shoulder.  With this choice, preprocessing the
  arm has provided the set $\{r_1,p_2\}$ and the third coordinate $a_3$ of
  $\namedvector{a}$. Notice that if $p_3=0$ then
  $\namedvector{a}=\coordinatevector{3}$, otherwise
  $\namedvector{a}=\vectorbycoordinates{1,0,r_1}$. Hence, if $a_3=1$ then
  $p_3=0$, but if $a_3>1$ then $p_3=1$ and $r_1=a_3$. In the $p_3=1$ case we
  get the integers $r_1$ and $p_2$, but in the case $p_3=0$, the integers
  $r_1$ and $p_2$ behave symmetrically, so we distinguish them arbitrarily.
  The algorithm finishes by postprocessing the arm.

  Notice that this algorithm covers not only the family $\ced_1$, but also
  some part of $\ced_2$. The remaining classes of $\ced_2$ will be discussed
  in~\ref{NEGY} (in accordance with this paragraph).
\end{para}

\subsection{The case $\ctr_1$.}\label{KETTO}

\begin{para} We assume that the arm of $v_1$ contains all vertices but one,
  which is not an ER-hand. Assume that the arm is in the direction $z_3$, and
  let $[(q_1,0,q_3),(0,p_2,p_3)]$ be its shoulder with $p_2\geq 2$, $q_1\geq 2$.
  Since $v_2$ is not an ER-vertex, $p_3>0$ and $q_3>0$
  (cf.~\ref{facts-hands}). If the third vertex of the face associated with
  $v_2$ is $(r_1,r_2,0)$, then $r_i>0$ ($i=1,2$) since otherwise $v_2$ would
  be in the arm of $v_1$. Therefore, $v_2$ corresponds to a central triangle
  with only crossing edges. Moreover, \eqref{eq:2} guarantees that $1\in
  \{r_1,p_3\}\cap \{r_2,q_3\}$. Let $(a_1,a_2,a_3)$ be the normal vector of the face
  of $v_2$.

  Let us collect some facts about such a Newton diagram in order to be able to
  find the right algorithm. Since $r_1=r_2=1$ is not possible (see
  \ref{alg:main-part}(\hyperlink{N0}{N=0})), we may assume that $p_3=1$. (This
  introduces a choice of the coordinates $z_1$ and $z_2$, and at this moment
  it is not clear how this choice fits with any property of $G^o$; this will
  be explained later.)

  We distinguish two cases.  The first case is $q_3=1$, then $v_2$ is a moving
  triangle, hence $a_1=p_2$ and $a_2=q_1$.  The second case is $q_3>1$, which
  we analyze in the rest of this paragraph.  Since $1\in\{r_2,q_3\}$, we get
  $r_2=1$. By \eqref{eq:19} one has
  \begin{equation}\label{NV}
    (a_1,a_2,a_3)=(p_2q_3-q_3+1,q_3r_1+q_1-r_1,r_1p_2+q_1-q_1p_2).
  \end{equation}
  From this and $q_3\geq 2$ one gets $a_1\geq 2p_2-1>p_2-1$ and $a_2\geq
  r_1+q_1>r_1,q_1$.  In particular, $a_1+a_2>p_2+q_1$ and hence
  \(\{a_1,a_2\}\neq\{p_1,q_2\}\).  The face value computed via the two vertices
  $(r_1,1,0)$ and $(0,p_2,1)$ gives $r_1a_1 + a_2 = p_2 a_2 + a_3$.
  Therefore, the integers \(r_1\), \(p_2\), \(q_1\), \(a_1\), \(a_2\) satisfy:
  \begin{equation}\label{CONST}
    \begin{gathered}
      a_3=r_1a_1-(p_2-1)a_2 \\
      0<r_1<a_2, \\
      0<p_2-1<a_1.
    \end{gathered}
  \end{equation}
\end{para}

\begin{para}[The algorithm] The two decorations of the legs of $v_2$ are
  $\calD=\{a_1,a_2\}$ (cf.~\eqref{eq:23}), where $\mygcd{a_1,a_2}=1$ by
  \ref{prop:2}(\ref{item:39}).  Preprocessing the arm (with hand $v_1$) has
  produced $a_3$ and the set $\calS=\{p_2,q_1\}$ (we cannot distinguish the two
  coordinates yet).  We shall compute the coordinates of \(v_2\) below, and
  then postprocess the arm to determine the rest of the Newton diagram.

  We distinguish two cases for computing \(v_2\).  First case: $\calD=\calS$.
  Let the two elements of this set be $a_1=p_2$ and $a_2=q_1$ (here is a
  choice between the \(z_1\) and \(z_2\) coordinates).  We select the
  d-minimal (as explained in~\ref{cor:str2}) solution $(r_1,r_2)$ of positive
  integers of the equation $r_1p_2+r_2q_1-q_1p_2=a_3$.  (The only reason for
  selecting the d-minimal solution is to obtain the d-minimal representative.)
  Then the vertices of $v_2$ are \((q_1,0,1)\), \((0,p_2,1)\),
  \((r_1,r_2,0)\).

  Second case: $\calD\neq\calS$.  We choose the unique \(6\)-tuple
  $(r_1,p_2,q_1,q_3,a_1,a_2)$ of positive integers with $\calS=\{a_1,a_2\}$,
  $\calD=\{p_2,q_1\}$ satisfying both \eqref{NV} and \eqref{CONST}.  (Uniqueness
  will be proved in the next paragraph.)  Then the vertices of $v_2$ (up to a
  permutation of the first two coordinates) are \((q_1,0,q_3)\), \((0,p_2,1)\)
  and \((r_1,1,0)\).
\end{para}
\begin{para}[Uniqueness of the \(6\)-tuple]
  Notice that once the choice between $a_1$ and $a_2$ is made, then
  \eqref{CONST} determines uniquely \(r_1\) and \(p_2\).  Then one gets
  \(q_1\) form \(\calD\) and also $q_3=(a_1-1)/(p_2-1)$.

  Assume for contradiction that by interchanging $a_1$ and $a_2$ we get
  another set of solutions \(\widetilde{r_1}\), \(\widetilde{p_2}\) and so on.
  Then, by \eqref{CONST}, $\widetilde{r_1}=a_1-p_2+1$ and
  $\widetilde{p_2}=a_2-r_1$.  Since $\widetilde{p_2}\in\{p_2,q_1\}$, there are two
  cases.

  If \(\widetilde{p_2} = q_1\) then substituting this in the expression of
  \(\widetilde{p_2}\) and using~\eqref{NV} for \(a_2\) produces
  \((q_3-2)r_1=-1\), whose left hand side is non-negative, a contradiction.

  On the other hand, if \(\widetilde{p_2}=p_2\) then from the expression of
  \(a_3\) in~\eqref{NV} (used for both sets of solutions) we obtain
  \(\widetilde{r_1}=r_1\).  Thus, again from~\eqref{NV}, we obtain \(a_1=a_2\)
  contradicting \(\mygcd{a_1,a_2}=1\).
\end{para}

\subsection{Two non-degenerate arms.}
\label{sec:two-non-degenerate}

\begin{para} Assume that there are either at least two vertices which are not
  in the arm of $v_1$, or there is only one such vertex, namely, \(v_2\).  In
  the latter case, we also assume that \(v_2\) is an ER-hand since the other
  case is treated in §\ref{KETTO}.  Anyway, $v_2$ is also a hand, so we
  preprocess its arm, too.  We face two cases: either the two arms (of $v_1$
  and $v_2$) cover all the vertices of $G^o$ (this fact characterizes the
  family $\ced_2$), or the arms contain all the vertices but one, which should
  be a central vertex/face (this is the family $\ctr_2$).
\end{para}

\begin{para}[The case $\ced_2$]\label{NEGY} If the arm of $v_2$
  contains all the vertices then we are in the situation of \ref{EGY}, and we
  are done. Assume that this is not the case. Fix the coordinates $z_i$ so
  that the arm of $v_i$ is in the direction of $z_i$ ($i=1,2$).  We select
  (arbitrarily) a common edge $\alpha=[(p,q,0),(0,0,c)]$ of the two arms, and let
  \(\bigtriangleup_i\) be the face adjacent to it in the direction \(z_i\).  In particular,
  \(\bigtriangleup_i\) lies in the arm of \(v_i\).  Let ${\namedvector{a}}^{(i)}$ be the
  normal vector of $\bigtriangleup_i$.  We seek the coordinates of these vectors and the
  edge \(\alpha\).

  By preprocessing the arms, we have obtained the sets $\{c,p\}$ and $\{c,q\}$,
  and the first two coordinates of both ${\namedvector{a}}^{(i)}$. By
  \eqref{eq:20} one has ${\namedvector{a}}^{(1)}×
  {\namedvector{a}}^{(2)}=(-p,-q,c)$, hence
  $c=a^{(1)}_1a^{(2)}_2-a^{(1)}_2a^{(2)}_1$. Hence we recover $\alpha$. Moreover,
  by face value computation, $a^{(i)}_3c=pa^{(i)}_1+qa^{(i)}_2$, hence we get
  the normal vectors as well. The algorithm finishes with postprocessing the
  arms.
\end{para}

\begin{para}[The case $\ctr_2$]\label{OT} Similarly as above, fix
  the coordinates $z_i$ so that the arm of $v_i$ is in the direction of $z_i$
  ($i=1,2$). We wish to determine the central triangle $\bigtriangleup$ using \ref{UJUJ},
  whose notations we will use.  Preprocessing the two non-degenerate arms, we
  have determined the sets $\{r_2,q_3\}$ and $\{p_3,r_1\}$, and the first two
  coordinates of $\normalvector{\bigtriangleup}$. The third coordinate of
  $\normalvector{\bigtriangleup}$ is the decoration of the leg adjacent to the vertex
  corresponding to $\bigtriangleup$, hence $\normalvector{\bigtriangleup}$ is known from $G^o$. For $1\leq
  i\leq 3$ denote by $\bigtriangledown_i$ the face of $\Gamma_+$ adjacent to $\bigtriangleup$ in the direction of
  the axis $z_i$.  Then, by the notations of \ref{MMM}, we already know the
  coordinates $a_1^{(1)}$ and $a_2^{(2)}$ from preprocessing the arms.
  Furthermore, $\bigtriangledown_3$ is a non-compact face with $a_3^{(3)}=0$ by~\ref{lem:10}.
  Therefore, \ref{MM} gives the face value $m_\bigtriangleup$. This, via the equations
  \eqref{eq:10} and \ref{cor:negdef} provide all the face values, in
  particular the face value of $\bigtriangledown_3$ too. This is $p_2q_1$. Since either $p_2$
  or $q_1$ is \(1\), we get the set $\{p_2,q_1\}$ as well. Hence, \ref{UJUJ}
  determines $\bigtriangleup$ (up to a permutation of coordinates). Then postprocessing the
  arms recovers the Newton diagram.
\end{para}

\section{Appendix}

\subsection{Some arithmetical properties of Newton boundaries.}\label{arith}

\begin{lemma}
  \label{lem:8}
  Let \(\bigtriangleup\) be a triangle whose vertices are lattice points.  Let
  \(\namedvector{a}\) and \(\namedvector{b}\) be the vectors of two of its
  sides.  Then
  \begin{equation}
    \label{eq:18}
    \namedvector{a} × \namedvector{b}=\pm
    \combinatorialarea{\bigtriangleup} \normalvector{\bigtriangleup}.
  \end{equation}
  In particular, if \(\bigtriangleup\) is an empty triangle with vertices \((0,p_2,p_3)\),
  \((q_1,0,q_3)\) and \((r_1,r_2,0)\), then
  \begin{equation}
    \label{eq:19}
    \normalvector{\bigtriangleup} = \vectorbycoordinates{p_2 q_3 + r_2 p_3 - r_2 q_3, q_3
      r_1 + p_3 q_1 - p_3 r_1, r_1 p_2 + q_1 r_2 - q_1 p_2}.
  \end{equation}
\end{lemma}

\begin{proof}
By the additivity of $\combinatorialarea{\bigtriangleup}$, we may assume that $\bigtriangleup$ is empty.
In that case \(\namedvector{a}\) and \(\namedvector{b}\) can be completed to a
base (see e.g.~\cite[p.~35]{Oda}), hence $\namedvector{a} × \namedvector{b}$
is primitive. The second part is a direct application. To verify the sign,
note that the scalar product of both vectors in \eqref{eq:19} with the
vertices of the triangle are positive.
\end{proof}

\begin{lemma}
  \label{lem:9}
  Let \(\bigtriangleup\) and \(\bigtriangledown\) be two adjacent lattice polygons.
  \begin{enumerate}[(a)]
  \item\label{item:77} Then the vector \(\namedvector{v}\) of their common
    edge is, up to a sign:
    \begin{equation}
      \label{eq:20}
      \namedvector{v}= \pm \frac{\multiplicity{\bigtriangleup}{\bigtriangledown}}{\innerdeterminant{\bigtriangleup}{\bigtriangledown}}
      \normalvector{\bigtriangleup} × \normalvector{\bigtriangledown}.
    \end{equation}
  \item\label{item:78} Assume that \(\bigtriangleup\) and \(\bigtriangledown\) are adjacent faces of a
    Newton polytope, \(\bigtriangleup\) is a triangle, and let \(\namedvector{a}\) be a
    vector from a point from their common edge to the third vertex of \(\bigtriangleup\).
    Then
    \begin{equation}
      \label{eq:21}
      \innerdeterminant{\bigtriangleup}{\bigtriangledown} =
      \frac{\multiplicity{\bigtriangleup}{\bigtriangledown}}{\combinatorialarea{\bigtriangleup}}
      \scalarproduct{\namedvector{a}}{\normalvector{\bigtriangledown}}.
    \end{equation}
  \item\label{item:79} Let \(\bigtriangleup\) be a triangle with vertices \((0,p_2,p_3)\),
    \((q_1,0,q_3)\) and \((q_1',0,q_3')\) with \(q_1' < q_1\), situated on a
    compact face of a Newton boundary.  Assume that \(\bigtriangleup\) has no lattice
    points other than its vertices and possible internal lattice points on its
    side on the \(z\sb{1} z\sb{3}\) plane.  Then the following expressions are
    equal and integers:
    \begin{equation}
      \label{eq:22}
      \frac{q_1-q_1'}{\multiplicity{\bigtriangleup}{\coordinatevector{2}}} =
      \frac{\scalarproduct{\normalvector{\bigtriangleup}}{\coordinatevector{3}}}{\innerdeterminant{\bigtriangleup}{\coordinatevector{2}}} \in \setN.
    \end{equation}
    In fact,
    \begin{equation}
      \label{eq:p2}
      \innerdeterminant{\bigtriangleup}{\coordinatevector{2}}=p_2.
    \end{equation}
  \end{enumerate}
\end{lemma}

\begin{proof} (\ref{item:77}) The vector \(\pm \namedvector{v}\) is
characterized by the fact that it is orthogonal to both normal vectors and it
is \(\multiplicity{\bigtriangleup}{\bigtriangledown}\) times a primitive vector.  The vector on the
right-hand side of \eqref{eq:20} has this property. For (\ref{item:78}), since
\(\namedvector{v}\) is orthogonal to \(\normalvector{\bigtriangledown}\),
Equations~\eqref{eq:20} and~\eqref{eq:18} give
\begin{equation*}
  \pm \combinatorialarea{\bigtriangleup}
  \frac{\innerdeterminant{\bigtriangleup}{\bigtriangledown}}{\multiplicity{\bigtriangleup}{\bigtriangledown}} \namedvector{v} =
  \pm \combinatorialarea{\bigtriangleup}\left(
    \normalvector{\bigtriangleup} ×
    \normalvector{\bigtriangledown} \right) =
  \left( \namedvector{a} × \namedvector{v}
  \right) × \normalvector{\bigtriangledown} =
  \scalarproduct{\namedvector{a}}{\normalvector{\bigtriangledown}} \namedvector{v}.
\end{equation*}
This gives \eqref{eq:21} up to a sign. Since scalar product of the normal
vector of a face assigns its minimum on the face (when restricted to the
Newton boundary), the scalar product in \eqref{eq:21} is positive, and hence
both sides of \eqref{eq:21} are positive.

For \eqref{eq:22}, we apply (\ref{item:77}) with
\(\normalvector{\bigtriangledown}=\coordinatevector{2}\) and \( \namedvector{v}=
\vectorbycoordinates{q_1-q_1',0,q_3-q_3'}\).  First notice that
\(\namedvector{v}\) is \(\multiplicity{\bigtriangleup}{\coordinatevector{2}}\) times a
primitive vector, hence $ (q_1-q_1')/\multiplicity{\bigtriangleup}{\coordinatevector{2}} =
\langle\namedvector{v}/ \multiplicity{\bigtriangleup}{\coordinatevector{2}},\coordinatevector{1}\rangle
\in \setN$.

On the other hand, taking scalar product of \eqref{eq:20} with
\(\coordinatevector{1}\), we obtain \eqref{eq:22} up to a sign.  Since both
expression are positive in \eqref{eq:22}, the sign is correct. The last
equality is a special case of (\ref{item:78}) with \(\bigtriangledown\) the \(z\sb{1}
z\sb{3}\) plane, because \(\combinatorialarea{\bigtriangleup} =
\multiplicity{\bigtriangleup}{\coordinatevector{2}}\).
\end{proof}

Recall (cf.~\ref{rem:uj}) that a non-compact face of $\Gamma$ (with \eqref{eq:2})
either lies on a coordinate plane, or it has an edge of type
\([(a,0,c),(0,1,b)]\) and normal vector \(\vectorbycoordinates{1,a,0}\) with
$a>0$.

\begin{lemma}
  \label{lem:10}
  Let an edge \(AB=[(a,0,c),(0,1,b)]\) lie on a compact face and on a
  non-compact one with normal vectors \(\namedvector{a}\) and
  \(\namedvector{n} \coloneqq \vectorbycoordinates{1,a,0}\), respectively
  ($a>0$).  Then
  \begin{equation}
    \label{eq:23}
    \innerdeterminant{\namedvector{n}}{\namedvector{a}} =
    \coordinateof{3}{\namedvector{a}}.
  \end{equation}
  Assume that \(C=(r,s,u)\) is a third vertex of the compact face, such that
  the triangle \(\trianglebyvertices{A}{B}{C}\) is empty. Then the determinant
  is also
  \begin{equation}
    \label{eq:24}
    \innerdeterminant{\namedvector{n}}{\namedvector{a}} =
    \scalarproduct{\vectorbycoordinates{r-a,s,u-c}}{\namedvector{n}} = r + (s-1)a.
  \end{equation}
\end{lemma}

\begin{proof}
Let \(\bigtriangleup\) be the empty triangle on the non-compact face with vertices:
\((a,0,c)\), \((0,1,b)\) and \((0,1,b+1)\) and \(\bigtriangledown\) denote the compact face.
Then \eqref{eq:21} with \(\namedvector{a}=\coordinatevector{3}\) yields
\eqref{eq:23}. The other equation is again an application of \eqref{eq:21}.
However, this time \(\bigtriangledown\) is the non-compact face, and \(\bigtriangleup\) is the triangle
with vertices \((a,0,c)\), \((0,1,b)\) and \((r,s,u)\).
\end{proof}

\begin{lemma}
  \label{lem:11}
  Let \(\bigtriangleup\) and \(\bigtriangledown\) be two adjacent triangular faces of a Newton diagram
  whose vertices lie on the coordinate planes containing the \(z\sb{3}\) axis.
  Further, let us assume that \(\bigtriangleup\) has an edge on the \(z\sb{1} z\sb{3}\)
  plane, which contains all the lattice points of the triangle except the
  third vertex. Let its determinant
  $\innerdeterminant{\bigtriangleup}{\coordinatevector{2}}$ be denoted by
  $\legdeterminant{\bigtriangleup}$.  Similarly, we suppose that \(\bigtriangledown\) has an edge \(\alpha\)
  either on the \(z\sb{1} z\sb{3}\) plane or on the \(z\sb{2} z\sb{3}\) plane
  containing all lattice points except the third vertex. Its determinant will
  be denoted by $\legdeterminant{\bigtriangledown}$. Then
  \begin{align}
    \label{eq:25}
    \alpha &\in \text{\(z\sb{1} z\sb{3}\) plane} &&\iff& \legdeterminant{\bigtriangleup} &=
    \legdeterminant{\bigtriangledown} \mid \innerdeterminant{\bigtriangleup}{\bigtriangledown},\\
    \label{eq:26}
    \alpha &\in \text{\(z\sb{2} z\sb{3}\) plane} &&\iff& \mygcd{\legdeterminant{\bigtriangleup},
      \legdeterminant{\bigtriangledown}, \innerdeterminant{\bigtriangleup}{\bigtriangledown}} &= 1.
  \end{align}
\end{lemma}
\begin{proof}
Let \(\namedvector{v}\) be the vector of the common edge of the triangles.
Let \(\namedvector{a}\) be the primitive vector parallel to the edge of \(\bigtriangleup\)
lying on the \(z\sb{1} z\sb{3}\) plane.  Finally, let \(\namedvector{c}\) be
the primitive vector parallel to \(\alpha\).  Now, \eqref{eq:20} combined with
\eqref{eq:18} implies that \(\innerdeterminant{\bigtriangleup}{\bigtriangledown}\) equals the triple
product \(\namedvector{a} \namedvector{v} \namedvector{c}\) (up to a sign).

If \(\alpha\) lies on the \(z\sb{1} z\sb{3}\) plane, \(\legdeterminant{\bigtriangleup} =
\legdeterminant{\bigtriangledown} = \coordinateof{2}{\namedvector{v}}\) by \eqref{eq:p2}.
Since the second coordinates of \(\namedvector{a}\) and \(\namedvector{c}\)
are \(0\), the number \(\coordinateof{2}{\namedvector{v}}\) divides the triple
product.  This proves the \(\implies\) part of~\eqref{eq:25}.

If \(\alpha\) lies on the \(z\sb{2} z\sb{3}\) plane, then
\(\coordinateof{2}{\namedvector{a}}=\coordinateof{1}{ \namedvector{c}}=0\).
Therefore, $\innerdeterminant{\bigtriangleup}{\bigtriangledown}$, modulo the greatest common divisor \(d\)
of \(\legdeterminant{\bigtriangleup} = \coordinateof{2}{\namedvector{v}}\) and
\(\legdeterminant{\bigtriangledown} = \coordinateof{1}{\namedvector{v}}\), is:
\begin{equation}
  \label{eq:27}
  \innerdeterminant{\bigtriangleup}{\bigtriangledown} =
  \namedvector{a} \namedvector{v} \namedvector{c} \equiv
  - \coordinateof{1}{\namedvector{a}} \cdot
  \coordinateof{3}{\namedvector{v}} \cdot
  \coordinateof{2}{\namedvector{c}} \pmod{d}.
\end{equation}
The three terms of the right-hand side are relative prime to \(d\) because
\(\normalvector{\bigtriangleup} = \pm \namedvector{a} × \namedvector{v}\), \(\normalvector{\bigtriangledown}
= \pm \namedvector{c} × \namedvector{v}\), and \(\namedvector{v}\) are
primitive. Hence the \(\implies\) part of~\eqref{eq:26} follows. We end the
proof by noticing that the right-hand sides of~\eqref{eq:25} and \eqref{eq:26}
are mutually exclusive.
\end{proof}

\subsection{The weighted homogeneous case (with one node).}\label{WHC}

Below $\mathcal{L}=(d_1,k_1;\dotsc;d_s,k_s)$ means that the unique vertex of
$G^o$ has $s$ leg-groups, the $i$th group has size $k_i\geq 1$ and decoration
$d_i>1$ (with $d_i\neq d_j$ for $i\neq j$, and $\sum_ik_i\geq 3$). The number $e$ is the
orbifold Euler number. One has the following cases:

\begin{enumerate}[1.]
\item $\mathcal{L}=(d,k)$

  Equation: $z_1^d+z_2^{k-1}z_3+z_2z_3^{k-1}$.

\item $\mathcal{L}=(d,2;D,2)$

  Equation: $z_1^dz_3+z_2^{2D}+z_3^2$, equivalently $z_1^{2d}+z_2^Dz_3+z_3^2$.

  (The equations are $\sim$-equivalent.)

\item $\mathcal{L}=(d,k;D,1)$, $d\mid D$

  Equation: $z_1^kz_2+z_1z_2^{(k-1)D/d+1}+z_3^d$.

\item $\mathcal{L}=(d,k;D,1)$, $\mygcd{d,D}=1$, $-edD=1$

  Equation: $z_1^d+ z_2^{(k-1)/D}z_3+ z_2z_3^k$.

\item $\mathcal{L}=(d,k;D,1)$, $\mygcd{d,D}=1$, $-edD=k$

  Equation: $z_1^dz_2+z_2^{(D+1)(k-1)/k}z_3+z_3^k$.

\item $\mathcal{L}=(a,2;b,2;c,2)$

  Equation: $z_1^{2a}+z_2^{2b}+z_3^{2c}$.

\item $\mathcal{L}=(a,k;b,1;c,1)$, $a\mid b$, $a\mid c$

  Equation: $z_1^{(bk)/a+1}z_2+z_1z_2^{(ck)/a+1}+z_3^a$.

\item $\mathcal{L}=(a,k;b,1;c,1)$, $b\mid c$ and $k>1$

  Equation: $z_1^az_2+z_2^{c/b+1}+z_3^{kb}$.

\item $\mathcal{L}=(a,k;b,1;c,1)$, $a\mid b$, $a\nmid c$, $-ebc=1$

  Equation: $z_1^{kc}z_2+z_2^{(bk)/a+1}+z_3^a$.

\item $\mathcal{L}=(a,k;b,1;c,1)$, $a\mid b$, $a\nmid c$, $A \coloneqq -ebc>1$

  Equation: $z_1^{(kc-1)/A+1}z_2+z_2^A+z_1z_3^a$.

\item $\mathcal{L}=(a,k;b,1;c,1)$, the numbers \(a\), \(b\), \(c\) does not
  divide each other, and $-eabc=k^2$

  Equation: $z_1^a+z_2^{kc}+z_3^{kb}$.

\item $\mathcal{L}=(a,k;b,1;c,1)$, the numbers \(a\), \(b\), \(c\) does not
  divide each other, $k=1$ and $A\coloneqq-eabc>1$

  Equation: $z_1^{(A-b)/a}z_2+z_2^{(A-c)/b}z_3+z_3^{(A-a)/c}z_1$, or\\*
  $z_1^{(A-b)/c}z_2+z_2^{(A-a)/b}z_3+z_3^{(A-c)/a}z_1$.

  (Only one of the equations have integer exponents, and this one gives the
  right diagram.)

\end{enumerate}

\bibliography{Newton-diagram}\bibliographystyle{amsplain}
\end{document}